\documentclass[a4paper, 10pt]{extarticle}
\usepackage{amsmath}
\usepackage{amssymb}
\usepackage{amsfonts}
\usepackage{graphicx}
\usepackage{subfigure}
\usepackage{authblk}
\usepackage{url}
\usepackage[textwidth=16cm,textheight=20cm]{geometry}
\usepackage{hyperref}
\hypersetup{
    colorlinks=true,
    citecolor=blue,
    linkcolor=blue,
    filecolor=blue,
    urlcolor=blue,
}
\newtheorem{theorem}{Theorem}

\newtheorem{example}[theorem]{Example}

\newtheorem{lemma}[theorem]{Lemma}

\newtheorem{remark}[theorem]{Remark}

\title{Differential Invariants in Algebra}
\author[1,\thanks{\textit{E-mail: }\texttt{valentin.lychagin@uit.no}}]{Valentin Lychagin}
\author[1\thanks{\textit{E-mail: }\texttt{mihail\underline{ }roop@mail.ru}}]{Mikhail Roop}
\affil[1]{V.A. Trapeznikov Institute of Control Sciences, Russian Academy of Sciences, 65 Profsoyuznaya Str., 117997 Moscow, Russia}

\begin{document}
\maketitle

\abstract{
In these lectures, we discuss two approaches to studying orbit spaces of algebraic Lie groups. Due to algebraic approach orbit space, or quotient, is an algebraic manifold, while from the differential viewpoint a quotient is a differential equation. The main goal of these lectures is to show that the differential approach gives us a better understanding of structure of invariants and orbit spaces. We illustrate this on classical equivalence problems, such as $\mathrm{SL}$ - classification of binary and ternary forms, and affine classification of algebraic plane curves.
}

\section{Introduction}
\label{sec:1}
The concept of an invariant appears whenever it comes to any kind of a classification problem. In these lectures, we would like to explain basic concepts of the invariant theory and show its applications to algebraic problems, such as $\mathrm{SL}$-classification of binary and ternary forms, and affine classification of algebraic plane curves. It seems helpful to us to recommend books \cite{VinPop, Alex} and references therein to the interested reader.

The origin of the invariant theory goes back to the middle of the 19th century and has not only mathematical motivation, such as affine classification of quadratic forms, finding canonical forms for equations of conics and quadrics, obtained in works of Euler, Lagrange, Cauchy, Gauss, but also a physical one (finding principal axes of inertia, investigation of planets' motion).

The first results on $\mathrm{SL}$-classification of binary forms belong to Boole (1841), who observed that discriminants of binary forms are invariant under linear transformations with determinant equal to 1. Later, in 1845, Cayley constructed invariants using the technique of hyperdeterminants developed by Cayley himself \cite{Cay-1,Cay-2}. In 1849, Aronhold provided a systematic study of ternary forms of degree 3, and two years later he gave a general formulation of invariant theory for algebraic forms. He also obtained differential equations for invariants of algebraic forms, that were also obtained by  Cayley for binary forms in 1852, which led to a series of works \cite{Cay-3,Cay-4,Cay-5,Cay-6} known as memoirs upon quantics.

In 1863, Aronhold observed that the number of rationally independent absolute invariants equals the difference between the number of coefficients of the form and the number of coefficients in a linear transformation (in modern terms, the difference between the dimension of the space of forms and the dimension of the group) \cite{Aron-1}. In 1861, Clebsch, using results of Aronhold, developed symbolic methods of finding invariants of algebraic forms \cite{Cleb-1}. These methods were later developed by Gordan and rapidly became popular.

In 1856, Cayley and Sylvester showed that binary forms of degrees up to four have a finite number of so-called \textit{irreducible covariants}. Covariant is a polynomial in $x$, $y$, and coefficients of the form, invariant under the transformations of the group (e.g. of $\mathrm{SL}_{2}$ transformations). Irreducibility means that such covariants cannot be expressed as rational functions of covariants of lower degree \cite{Sil-1}. This became the origin of the finiteness problem for generating set of invariants.

Gordan was the first who proved the finiteness of a number of covariants for the binary form of arbitrary degree (Gordan's theorem) \cite{Gord}, and his method allowed to construct a complete system of irreducible covariants for binary forms of degrees 5 and 6. Later, Sylvester discovered the same result for the case of a binary form of degree 12. In 1880, von Gall constructed a complete system of covariants for a binary form of degree 8, and eight years later for that of degree 7, which turned out to be more complicated than the case of degree 8 \cite{Gall-1,Gall-2}. Binary forms of degree 7 were also elaborated by Dixmier and Lazard \cite{Dix}. Hammond provided the proof for the case of binary seventhics \cite{Ham}.

Finally, in 1890, Hilbert gave a complete proof of Gordan's result for the case of arbitrary $n$-ary forms of an arbitrary degree \cite{Hilb-1}.

While solving the problem of constructing a complete system of irreducible invariants and covariants, the very notion of an \textit{invariant} was changing. The theory of \textit{differential invariants} was developed by Halphen in 1878 in his thesis \cite{Halph} and was later generalized by Norwegian mathematician Sophus Lie, who showed that all previous results of invariant theory are particular cases of more general theory of invariants of continuous transformation groups \cite{Lie-1,Lie-2}. Lie did not use symbolic methods of Aronhold and Clebsch, that hardly could be extended to the cases of binary forms of higher degrees due to their dramatic bulkiness.

In the context of modern invariant theory and simultaneously in the context of these lectures, it is worth mentioning such results as Rosenlicht \cite{Ros} and global Lie-Tresse theorems \cite{KrugLych}, that justified the appearance of rational differential invariants in classification problems and paved a way for solving algebraic equivalence problems using differential-geometric techniques \cite{LychBib-1,LychBib-2}. This will be the core point of the present lectures.

The paper is organized as follows. In Sect. \ref{sec:2}, we start with $\mathrm{SL}_{2}(\mathbb{C})$ classification of binary forms and explain how to get rational differential invariants using the observation that binary forms are solutions of the Euler equation. In Sect. \ref{sec:3}, we give a general introduction to modern invariant theory together with discussion of Rosenlicht and Lie-Tresse theorems and explanation how the last can be used to find smooth solutions to PDEs, as well as those with singularities. Sect. \ref{sec:4} is devoted to affine classification of algebraic plane curves. The last Sect. \ref{sec:5} concerns the problem of $\mathrm{SL}_{3}(\mathbb{C})$-classification of ternary forms using results obtained in the previous sections. 

All essential computations for this paper were performed in Maple with the DifferentialGeometry package created by I. Anderson and his team \cite{AndTor}, and the first author is grateful to him for the very first introduction to the package.
\section{Invariants of Binary Forms}
\label{sec:2}
In this section, we study $\mathrm{SL}_{2}$ - invariants of binary $n$ - forms. We show the difference between algebraic and differential approaches and the power of differential one in finding invariants.
\subsection{Algebraic Point of View}
Binary form of degree $n$ is a homogeneous polynomial on $\mathbb{C}^{2}$
\begin{equation}
\label{bin-f}
\phi_{b}=\sum\limits_{i=0}^{n}b_{i,n-i}\frac{x^{i}}{i!}\frac{y^{n-i}}{(n-i)!},\quad b_{i,n-i}\in\mathbb{C}.
\end{equation}
The space of all binary forms of degree $n$ is $\mathcal{B}_{n}\simeq\mathbb{C}^{n+1}$. The action of the Lie group
\begin{equation*}
\mathrm{SL}_{2}(\mathbb{C})=\left\{A\in\mathrm{Mat}_{2\times 2}(\mathbb{C})\mid \det(A)=1\right\}
\end{equation*}
on $\mathcal{B}_{n}$ is defined by the following way:
\begin{equation}
\label{SL2-act}
A\colon\mathcal{B}_{n}\ni\phi_{b}\mapsto A\phi_{b}=\phi_{b}\circ A^{-1}\in\mathcal{B}_{n}.
\end{equation}
This action induces the action on coefficients $b_{i,n-i}$. Due to algebraic approach, where we believe that the quotient is an algebraic manifold, to describe the quotient space $\mathcal{B}_{n}/\mathrm{SL}_{2}(\mathbb{C})$ one needs to find polynomials $I(b)=I(b_{0,n},\ldots,b_{n,0})$ invariant under the action \eqref{SL2-act}. Such functions are called \textit{algebraic invariants}.
\begin{theorem}[Gordan-Hilbert, \cite{Gord,Hilb-1}]
The algebra of polynomial $\mathrm{SL}_{2}$ - invariants of binary $n$-forms is finitely generated, and the quotient space is an affine, algebraic manifold.
\end{theorem}
However, the problem of finding generators of this algebra and syzygies in this algebra turned out to be specific for every $n$. For instance, the case of $n=3$ was elaborated by Bool in 1841, who observed that the discriminant of the cubic is an invariant. This became the origin of the classical invariant theory. Results regarding the case of $n=4$ belong to Bool, Cayley and Eisinsteine (1840-1850) \cite{Cay-1,Cay-2,Hilb,Shur}. For quintic $(n=5)$, the invariants were found by Sylvester and Hilbert (see, for example, \cite{Hilb,Shur}). They are dramatically huge to write down explicitly, the invariant of degree 18 found by Hermite contains 848 terms! The main problem is that there is no general approach in the classical invariant theory. This motivates us to develop a differential approach \cite{LychBib-1,LychBib-2}.
\subsection{Differential Point of View}
The key idea underlying the differential approach is to identify $\mathcal{B}_{n}$ with the space of smooth solutions to Euler equation
\begin{equation}
\label{euler}
xf_{x}+yf_{y}=nf.
\end{equation}
It is worth mentioning that class of solutions to \eqref{euler} includes not only binary $n$-forms, but also other homogeneous functions of degree $n$. Thus, solving the problem for all solutions to \eqref{euler} we at the same time solve the problem of $\mathrm{SL}_{2}$-equivalence of binary forms.

Equation \eqref{euler} defines a smooth submanifold $\mathcal{E}_{1}$ in the space of 1-jets $\mathbf{J}^{1}=J^{1}\left(\mathbb{C}^{2}\right)$ of functions on $\mathbb{C}^{2}$:
\begin{equation*}
\mathcal{E}_{1}=\left\{xu_{10}+yu_{01}=nu_{00}\right\}\subset\mathbf{J}^{1}.
\end{equation*}
Solutions of \eqref{euler} are special type surfaces $L_{f}\subset\mathcal{E}_{1}$
\begin{equation*}
L_{f}=\left\{u_{00}=f(x,y),\,u_{10}=f_{x},\, u_{01}=f_{y}\right\}\subset\mathcal{E}_{1}.
\end{equation*}
It is often reasonable to consider not only equation \eqref{euler}, but also a collection of its differential consequences up to some order $k$, i.e. a prolongation $\mathcal{E}_{k}\subset\mathbf{J}^{k}$. The space $\mathbf{J}^{k}$ is a space of $k$-jets of smooth functions on $\mathbb{C}^{2}$:
\begin{equation*}
\mathbf{J}^{k}=\left\{[f]_{p}^{k}\mid p\in\mathbb{C}^{2},\,f\in C^{\infty}\left(\mathbb{C}^{2}\right)\right\},
\end{equation*}
where $[f]_{p}^{k}$ is the equivalence class of functions, whose Taylor polynomials of the length $k$ at the point $p\in\mathbb{C}^{2}$ are the same (values and all derivatives up to order $k$ at the point $p$ coincide). The space of $k$-jets is equipped with canonical coordinates $(x,y,u_{00},\ldots,u_{ij},\ldots)$, $0\le i+j\le k$, $\dim\left(\mathbf{J}^{k}\right)=\binom{k+2}{2}+2$, and
\begin{equation*}
u_{ij}\left([f]_{p}^{k}\right)=\frac{\partial^{i+j}f}{\partial x^{i}\partial y^{j}}(p).
\end{equation*}

The action $A\colon\mathbb{C}^{2}\to\mathbb{C}^{2}$ of the group $\mathrm{SL}_{2}$ can be prolonged to $\mathbf{J}^{k}$ by the natural way
\begin{equation*}
A^{(k)}\colon\mathbf{J}^{k}\to\mathbf{J}^{k},\quad A^{(k)}\left([f]_{p}^{k}\right)=[Af]_{Ap}^{k}.
\end{equation*}
Moreover, if
\begin{equation*}
L_{f}^{(k)}=\left\{u_{ij}=\frac{\partial^{i+j}f}{\partial x^{i}\partial y^{j}},\, 0\le i+j\le k\right\}
\end{equation*}
is a graph of the $k$-jet of function $f$, then
\begin{equation*}
A^{(k)}\left(L_{f}^{(k)}\right)=L_{Af}^{(k)}.
\end{equation*}
Let us now put $k=n$ and let $\mathcal{E}_{n}\subset\mathbf{J}^{n}$ be the $(n-1)$-prolongation of the Euler equation together with $u_{ij}=0$:
\begin{equation*}
\mathcal{E}_{n}=\left\{\frac{d^{k+l}}{dx^{k}dy^{l}}\left(xu_{10}+yu_{01}-nu_{00}\right)=0,\,0\le k+l\le n-1,\, u_{ij}=0,\, n+1\le i+j\right\}.
\end{equation*}
One can show that $\dim\mathcal{E}_{n}=n+3$. The prolongations $A^{(n)}$ of group elements $A\in\mathrm{SL}_{2}$ preserve the submanifold $\mathcal{E}_{n}$ and therefore define the action $A^{(n)}\colon\mathcal{E}_{n}\to\mathcal{E}_{n}$. Since $L_{\phi}^{(n)}\subset\mathcal{E}_{n}$, any binary $n$-form can be considered as a solution to $\mathcal{E}_{n}$. The property $A^{(n)}\left(L_{\phi}^{(n)}\right)=L_{A\phi}^{(n)}$ shows that the group $\mathrm{SL}_{2}(\mathbb{C}^{2})$ is a symmetry group of the Euler equation.

A rational function $I\in C^{\infty}(\mathcal{E}^{k})$ is said to be a \textit{rational differential $\mathrm{SL}_{2}$-invariant of order $k$}, or simply \textit{differential invariant}, if $I\circ A^{(k)}=I$, for all $A\in\mathrm{SL}_{2}(\mathbb{C})$.

As we shall see further, the Lie-Tresse theorem states that the algebra of rational differential $\mathrm{SL}_{2}$-invariants of order $\le n$ on the Euler equation $\mathcal{E}_{n}$ gives us realization of the quotient $\mathcal{E}_{n}/\mathrm{SL}_{2}(\mathbb{C})$ as a new differential equation of order 3, and $\mathrm{SL}_{2}(\mathbb{C})$-orbits of binary $n$-forms correspond to solutions of this equation.

The following observations will be important for us.
\begin{itemize}
\item the plane $\mathbb{C}^{2}$ is the affine space, i.e. a space with the standard translation of vectors (trivial connection) and distinguished point $\mathbf{0}$
\item the plane $\mathbb{C}^{2}$ is the symplectic space, equipped with the structure form $\Omega=dx\wedge dy$
\item the group $\mathrm{SL}_{2}(\mathbb{C})$ preserves these both affine and symplectic structures, and the point $\mathbf{0}$.
\end{itemize}
As we shall see further, these structures will allow us to equip the set of differential $\mathrm{SL}_{2}(\mathbb{C})$-invariants with additional structures and will give us explicit methods of finding invariants.
\subsection{Relations between Algebraic and Differential Invariants}
One can easily see that due to \eqref{bin-f}
\begin{equation*}
b_{i,n-i}=\frac{\partial^{n}\phi_{b}}{\partial x^{i}\partial y^{n-i}}.
\end{equation*}
Therefore, the function $I(b_{n,0},\ldots,b_{0,n})$ is an $\mathrm{SL}_{2}(\mathbb{C})$-invariant if and only if $I(u_{n0},\ldots,u_{0n})$ is a differential $\mathrm{SL}_{2}(\mathbb{C})$-invariant of order $n$. Thus, algebraic $\mathrm{SL}_{2}(\mathbb{C})$-invariants of binary $n$-forms are differential invariants of the form $I(u_{0n},\ldots,u_{n0})$ and finding differential invariants we simultaneously find also algebraic ones.
\subsection{Lie Equation}
Since the Lie group $\mathrm{SL}_{2}(\mathbb{}C)$ is connected, the condition $I\circ A^{(k)}=I$ can be written in an infinitesimal form:
\begin{equation}
\label{lie-eq}
X^{(k)}(I)=0,\quad X\in\mathfrak{sl}_{2},
\end{equation}
where $X^{(k)}$ is the $k$th prolongation of the vector field $X\in\mathfrak{sl}_{2}$, and equation \eqref{lie-eq} is called \textit{Lie equation}. The Lie algebra $\mathfrak{sl}_{2}$ is generated by vector fields
\begin{equation*}
\mathfrak{sl}_{2}=\langle X_{+}=x\partial_{y},\,X_{-}=y\partial_{x},\,X_{0}=x\partial_{x}-y\partial_{y}\rangle
\end{equation*}
with commutators
\begin{equation}
\label{alg-str}
[X_{+},X_{-}]=X_{0},\quad [X_{0},X_{+}]=2X_{+},\quad [X_{0},X_{-}]=-2X_{-}.
\end{equation}
Due to Lie algebra structure \eqref{alg-str}, condition $X_{0}^{(k)}(I)=0$ is not independent, and Lie equation \eqref{lie-eq} becomes
\begin{equation*}
X_{+}^{(k)}(I)=0,\quad X_{-}^{(k)}(I)=0.
\end{equation*}
This equation also appeared in Hilbert's lectures \cite{Hilb}.

Following some empirical observations, according to which the number of functionally independent invariants equals the codimension of the regular orbit (we shall explain this strictly by means of the Rosenlicht theorem in the forthcoming sections), let us now compute the numbers of functionally independent algebraic and differential invariants.

Since
\begin{equation*}
\dim(\mathbf{J}^{k})=\frac{(k+1)(k+2)}{2}+2,
\end{equation*}
the number of independent differential invariants of $k$th order on $\mathbf{J}^{k}$ equals
\begin{equation*}
\dim(\mathbf{J}^{k})-\dim(\mathfrak{sl}_{2})=\frac{k(k+3)}{2}.
\end{equation*}
Since $\dim(\mathcal{E}_{n})=n+3$, the number of differential invariants of binary $n$-forms equals $\dim(\mathcal{E}_{n})-3=n$, and the number of independent algebraic invariants of binary $n$-forms equals $\dim(\mathbb{C}^{n+1})-3=n+1-3=n-2$.

This discussion is true for the case $n\ge 3$, when the Lie algebra of the stabilizer of the form is trivial. In the case $n=2$ its dimension equals 1, and therefore there is only one invariant in this case, which is the discriminant.

\subsection{Resultants and Discriminants}
Here, we will repeat the Boole's result on the $\mathrm{SL}_{2}$-invariance of the discriminant of binary forms.

Any binary $n$-form can be represented as a product of linear functions $I_{i}^{\phi}$, $i=1,\ldots,n$:
\begin{equation*}
\phi=\prod\limits_{i=1}^{n}I_{i}^{\phi}.
\end{equation*}
Obviously, functions $I_{i}^{\phi}$ are defined up to multipliers $\lambda_{i}$: $I_{i}^{\phi}\mapsto\lambda_{i}I_{i}^{\phi}$, where $\prod\limits_{i=1}^{n}\lambda_{i}=1$. Let $\psi\in\mathcal{B}_{n}$ be another binary form, $\psi=\prod\limits_{i=1}^{m}I_{i}^{\psi}$. Then, one can define \textit{resultant} between forms $\phi$ and $\psi$ by the following way:
\begin{equation*}
\mathrm{Res}(\phi,\psi)=\prod\limits_{i,j}[I_{i}^{\phi},I_{j}^{\psi}],
\end{equation*}
where $[I_{i}^{\phi},I_{j}^{\psi}]$ is the Poisson bracket associated with the symplectic form $\Omega=dx\wedge dy$.

The function
\begin{equation*}
\mathrm{Discr}(\phi)=\mathrm{Res}(\phi_{x},\phi_{y}),
\end{equation*}
is called \textit{discriminant}.

Remark that here $(x,y)$ are canonical coordinates of the vector space $\mathbb{C}^{2}$, i.e. $\Omega=dx\wedge dy$ in these coordinates.

Let us collect basic properties of discriminants and resultants.
\begin{enumerate}
\item $\mathrm{Res}(\phi,\psi)$ does not depend on scalings $I_{i}^{\phi}\mapsto\alpha_{i}I_{i}^{\phi}$, $I_{i}^{\psi}\mapsto\beta_{i}I_{i}^{\psi}$
\item $\mathrm{Res}(\phi,\psi)$ is a polynomial in coefficients of $\phi$, $\psi$ of degree $(n+m)$
\item $\mathrm{Res}(\phi,\psi)$ is an $\mathrm{SL}_{2}(\mathbb{C})$-invariant: $\mathrm{Res}(A\phi,A\psi)=\mathrm{Res}(\phi,\psi)$
\item $\mathrm{Discr}(\phi)$ is a polynomial $\mathrm{SL}_{2}(\mathbb{C})$-invariant of degree $(2n-2)$.
\end{enumerate}
Using discriminants and resultants one gets algebraic invariants from differential ones.
\begin{example}
Consider the following binary form of degree 3:
\begin{equation}
\label{cubic-ex}
\phi_{3}(x,y)=x^{3}+a_{1}x^{2}y+a_{2}xy^{2}+a_{3}y^{3}
\end{equation}
\begin{enumerate}
\item The discriminant $\mathrm{Discr}(\phi)$ of cubic \eqref{cubic-ex}
\begin{equation*}
J_{1}=\mathrm{Discr}(\phi)=12a_{1}^{3} a_{3}-3 a_{1}^{2} a_{2}^{2}-54 a_{1} a_{2} a_{3}+12 a_{2}^{3}+81 a_{3}^{3}
\end{equation*}
is a polynomial $\mathrm{SL}_{2}(\mathbb{C})$-invariant of order 4. This illustrates the property 4.

\item Let us take the differential $\mathrm{SL}_{2}$-invariant $u_{20}u_{02}-u_{11}^{2}$ and restrict it on the cubic \eqref{cubic-ex}. We get the following quadric
\begin{equation*}
\phi_{2}(x,y)=4(3 a_{2}- a_{1}^{2})x^{2}+4(9 a_{3}- a_{1} a_{2})xy+4(3 a_{1} a_{3}- a_{2}^{2})y^{2}.
\end{equation*}
Taking its discriminant, we get the polynomial invariant $J_{2}=-16J_{1}$. This illustrates how one can get polynomial invariants from differential ones.
\end{enumerate}
\end{example}
\subsection{Operations and Structures on Invariants}
\subsubsection{Monoid Structure}
Any function $\phi\in C^{\infty}(\mathbf{J}^{k})$ generates a differential operator by the following way:
\begin{equation*}
\widehat{\phi}\colon C^{\infty}(\mathbb{C}^{2})\to C^{\infty}(\mathbb{C}^{2}),\quad
\end{equation*}
or in coordinates
\begin{equation*}
\widehat{\phi}\colon f(x,y)\mapsto\phi\left(x,y,f,f_{x},f_{y},\ldots\right),
\end{equation*}
if $\phi=\phi(x,y,u_{00},u_{10},u_{01},\ldots)$. Then, condition for $\phi$ to be an $\mathrm{SL}_{2}(\mathbb{C})$-invariant reads
\begin{equation*}
A\circ\widehat{\phi}=\widehat{\phi}\circ A, \quad A\in\mathrm{SL}_{2}(\mathbb{C}).
\end{equation*}
Now we can introduce an operation $\ast$ of composition for invariants by the following way:
\begin{equation*}
\widehat{\phi\ast\psi}=\widehat{\phi}\circ\widehat{\psi}.
\end{equation*}
\begin{example}
\begin{equation*}
u_{00}\ast\psi=\psi,\quad u_{10}\ast\psi=\frac{d\psi}{dx},\quad u_{01}\ast\psi=\frac{d\psi}{dy},\quad u_{ij}\ast\psi=\frac{d^{i+j}\psi}{dx^{i}dy^{j}},
\end{equation*}
\begin{equation*}
(u_{20}u_{02}-u_{11}^{2})\ast\psi=\frac{d^{2}\psi}{dx^{2}}\frac{d^{2}\psi}{dy^{2}}-\left(\frac{d^{2}\psi}{dxdy}\right)^{2},
\end{equation*}
where
\begin{equation*}
\frac{d}{dx}=\frac{\partial}{\partial x}+\sum\limits_{i,j=0}u_{i+1,j}\frac{\partial}{\partial u_{ij}},\quad \frac{d}{dy}=\frac{\partial}{\partial y}+\sum\limits_{i,j=0}u_{i,j+1}\frac{\partial}{\partial u_{ij}}
\end{equation*}
are total derivatives.
\end{example}
Note that the composition of differential invariants of orders $k$ and $l$ is a differential invariant of order $(k+l)$, and composition with $u_{00}$ gives us the same invariant. This means that the composition operation endows the set of differential $\mathrm{SL}_{2}(\mathbb{C})$-invariants with a monoid structure.
\begin{theorem}
The set of differential $\mathrm{SL}_{2}(\mathbb{C})$-invariants is a monoid with unit $u_{00}$.
\end{theorem}
\begin{example}
The differential $\mathrm{SL}_{2}(\mathbb{C})$-invariants of order 1 are
\begin{equation*}
\phi=F(u_{00},xu_{10}+yu_{01}).
\end{equation*}
Let $\psi$ be another invariant of order $k$. Then,
\begin{equation*}
\phi\ast\psi=F\left(\psi,x\frac{d\psi}{dx}+y\frac{d\psi}{dy}\right)
\end{equation*}
is a differential invariant of order $(k+1)$.
\end{example}
\subsubsection{Poisson Structure}
Recall that the symplectic form $\Omega=dx\wedge dy$ is $\mathrm{SL}_{2}$-invariant. Define the Poisson bracket for functions on jet spaces by the following way:
\begin{equation*}
\widehat{d}\phi\wedge \widehat{d}\psi=[\phi,\psi]\Omega,
\end{equation*}
where $\widehat{d}f=\frac{df}{dx}dx+\frac{df}{dy}dy$ is the total differential, $f\in C^{\infty}(\mathbf{J}^{k})$. As we shall see below, $\widehat{d}$ is an invariant operator. Then, we get
\begin{equation*}
[\phi,\psi]=\frac{d\phi}{dx}\frac{d\psi}{dy}-\frac{d\phi}{dy}\frac{d\psi}{dx},
\end{equation*}
and if $\phi$ and $\psi$ are differential $\mathrm{SL}_{2}$-invariants, then $[\phi,\psi]$ is a differential invariant too.
\begin{theorem}
The algebra of $\mathrm{SL}_{2}$-invariants is a Poisson algebra.
\end{theorem}
\begin{example}
Let us take two differential $\mathrm{SL}_{2}(\mathbb{C})$-invariants: $J_{1}=u_{00}$ and $J_{2}=u_{20}u_{02}-u_{11}^{2}$. Taking the Poisson bracket between them we get a differential $\mathrm{SL}_{2}(\mathbb{C})$-invariant of the third order:
\begin{equation*}
J_{3}=[J_{1},J_{2}]=u_{01}(2u_{11}u_{21}-u_{02}u_{30}-u_{20}u_{12})+u_{10}(u_{02}u_{21}+u_{20}u_{03}-2u_{11}u_{12}).
\end{equation*}
As en exercise, we propose to check it to the reader.
\end{example}
\subsubsection{Invariant Frame}
Taking the $k$th term in the Taylor decomposition of a function $f(x,y)$, we get symmetric differential forms
\begin{equation*}
d_{k}f=\sum\limits_{i=0}^{k}\frac{\partial^{k}f}{\partial x^{i}\partial y^{k-i}}\frac{dx^{i}}{i!}\frac{dy^{k-i}}{(k-i)!},\quad k=1,2,\ldots
\end{equation*}
We shall see later on that these tensors are defined by the affine connection, which is in our case the trivial connection. Therefore, they are invariants of the affine transformations, i.e.
\begin{equation*}
d_{k}(Af)=A(d_{k}f),\quad A\in\mathrm{SL}_{2}(\mathbb{C}).
\end{equation*}
Let us define tensors $\Theta_{k}$ on jet spaces by the following way:
\begin{equation*}
\Theta_{k}=\sum\limits_{i=0}^{k}u_{i,k-i}\frac{dx^{i}}{i!}\frac{dy^{k-i}}{(k-i)!}.
\end{equation*}
Then, $d_{k}f=\Theta_{k}|_{L_{f}^{k}}$, and $\Theta_{k}$ are $\mathrm{SL}_{2}$-invariants.

On the space $\mathbf{J}^{2}$ we have the following $\mathrm{SL}_{2}$-invariant tensors:
\begin{eqnarray*}
  \Theta_{1} &=& u_{10}dx+u_{01}dy, \\
  \Theta_{2} &=& u_{20}\frac{dx^2}{2}+u_{11}dxdy+u_{02}\frac{dy^2}{2},\\
  \Omega&=&dx\wedge dy.
\end{eqnarray*}
As we shall see further, the Lie-Tresse theorem states that the algebra of differential invariants is a differential algebra, and we now turn the algebra of invariants into the differential algebra by introducing the invariant derivations
\begin{equation*}
\nabla_{i}=A_{i}\frac{d}{dx}+B_{i}\frac{d}{dy},\quad i=1,2,
\end{equation*}
where $A_{i}$ and $B_{i}$ are functions on $\mathbf{J}^{2}$, satisfying the conditions:
\begin{equation*}
\nabla_{1}\rfloor\Omega=\Theta_{1},\quad\nabla_{2}\rfloor\Theta_{2}=\Theta_{1}.
\end{equation*}
Direct computations give us the following result:
\begin{eqnarray}
\label{nablas}
  \nabla_{1} &=& u_{01}\frac{d}{dx}-u_{10}\frac{d}{dy}, \\
  \nabla_{2} &=& \frac{2(u_{02}u_{10}-u_{11}u_{01})}{\Delta_{2}}\frac{d}{dx}+\frac{2(u_{20}u_{01}-u_{11}u_{10})}{\Delta_{2}}\frac{d}{dy},
\end{eqnarray}
where $\Delta_{2}=u_{20}u_{02}-u_{11}^{2}$.

Their bracket is
\begin{equation*}
[\nabla_{1},\nabla_{2}]=A\nabla_{1}+B\nabla_{2},
\end{equation*}
where $A$ and $B$ are differential $\mathrm{SL}_{2}$-invariants of order 3, and
\begin{equation*}
A|_{\mathcal{E}_{3}}=\frac{2(2-n)}{n-1},\quad B|_{\mathcal{E}_{3}}=0.
\end{equation*}
\begin{theorem}
Let $\phi$ be a differential $\mathrm{SL}_{2}$-invariant of order $\le k$. Then, $\nabla_{1}(\phi)$ and $\nabla_{2}(\phi)$ are differential $\mathrm{SL}_{2}$-invariants of order $\le k+1$.
\end{theorem}
This means that the algebra of differential $\mathrm{SL}_{2}$-invariants equipped with invariant derivations $\nabla_{1}$ and $\nabla_{2}$ becomes a differential algebra. Summarizing all above discussion, we have:
\begin{theorem}
The algebra of differential $\mathrm{SL}_{2}$-invariants is a
\begin{itemize}
\item monoid with unit $u_{00}$
\item Poisson algebra
\item differential algebra
\end{itemize}
\end{theorem}
We can see that the differential viewpoint allows us to endow the set of invariants with much more interesting structures comparing with those we had in the algebraic situation.
\subsection{Invariant coframe}
Let us now construct the dual frame $\langle\omega_{1},\omega_{2}\rangle$, which is an $\mathrm{SL}_{2}$-invariant coframe, where $\omega_{i}=a_{i}dx+b_{i}dy$ and coefficients $a_{i}$, $b_{i}$ are such that $\omega_{i}(\nabla_{j})=\delta_{ij}$.

Simple computations give us
\begin{eqnarray*}
  \omega_{1} &=& \frac{u_{20}u_{01}-u_{11}u_{10}}{J_{21}}dx-\frac{u_{02}u_{10}-u_{11}u_{01}}{J_{21}}dy, \\
  \omega_{2} &=& \frac{\Delta_{2}}{2J_{21}}(u_{10}dx+u_{01}dy),
\end{eqnarray*}
where
\begin{equation*}
J_{21}=u_{01}^{2}u_{20}-2u_{10}u_{01}u_{11}+u_{10}^{2}u_{02}
\end{equation*}
is an $\mathrm{SL}_{2}$-invariant of order 2, called \textit{flex invariant} \cite{LychGold}.

The original coframe $\langle dx,dy\rangle$ is expressed in terms of $\langle\omega_{1},\omega_{2}\rangle$ as
\begin{eqnarray*}
 dx &=&u_{01}\omega_{1}+\frac{2(u_{02}u_{10}-u_{11}u_{01})}{\Delta_{2}}\omega_{2}, \\
 dy &=&-u_{10}\omega_{1}+\frac{2(u_{20}u_{01}-u_{11}u_{10})}{\Delta_{2}}\omega_{2}.
\end{eqnarray*}
And finally we are able to write down the invariant tensors $\Theta_{k}$ in the form
\begin{equation*}
\Theta_{k}=\sum\limits_{i=0}^{k}I_{i,k-i}\frac{\omega_{1}^{i}\omega_{2}^{k-i}}{i!(k-i)!}.
\end{equation*}
Since $\Theta_{k}$ are invariants, $\omega_{1,2}$ are invariants, we get:
\begin{theorem}
Functions $I_{i,j}$ are $\mathrm{SL}_{2}$-invariants of order $(i+j)$, and any rational differential invariant is a rational function of them.
\end{theorem}
\begin{example}
\begin{itemize}
\item $k=0$

The only invariant of the zeroth order is $I_{0,0}=u_{00}$.
\item $k=1$
\begin{equation*}
\Theta_{1}=\frac{2J_{21}}{\Delta_{2}}\omega_{2}.
\end{equation*}
\item $k=2$
\begin{equation*}
\Theta_{2}=\frac{J_{21}}{2}\omega_{1}^{2}+\frac{2J_{21}}{\Delta_{2}}\omega_{2}^{2}.
\end{equation*}
\item $k=3$
\begin{equation*}
I_{3,0}=-\frac{1}{6}u_{0 3}u_{1 0}^3+\frac{1}{2}u_{1 2}u_{0 1}u_{1 0}^2-\frac{1}{2}u_{2 1}u_{0 1}^2u_{1 0}+\frac{1}{6}u_{0 1}^3u_{3 0},
\end{equation*}
\begin{equation*}
\begin{split}
I_{1,2}&=\Delta_{2}^{-2}((2u_{1 1}^2u_{3 0}-4u_{1 1}u_{2 0}u_{2 1}+2u_{1 2}u_{2 0}^2)u_{0 1}^3+2u_{1 0}(u_{2 1}u_{1 1}^2-{}\\&-2u_{0 2}u_{3 0}u_{1 1}+u_{2 0}(2u_{2 1}u_{0 2}-u_{0 3}u_{2 0}))u_{0 1}^2+2u_{1 0}^2(u_{0 2}^2u_{3 0}-{}\\&-2u_{0 2}u_{1 2}u_{2 0}+2u_{0 3}u_{1 1}u_{2 0}-u_{1 1}^2u_{1 2})u_{0 1}-2u_{1 0}^3(u_{0 2}^2u_{2 1}-2u_{0 2}u_{1 1}u_{1 2}+u_{0 3}u_{1 1}^2)),
\end{split}
\end{equation*}
\begin{equation*}
\begin{split}
I_{2,1}&=\Delta_{2}^{-1}((-u_{1 1}u_{3 0}+u_{2 0}u_{2 1})u_{0 1}^3+u_{1 0}(u_{0 2}u_{3 0}+u_{1 1}u_{2 1}-2u_{1 2}u_{2 0})u_{0 1}^2-{}\\&-u_{1 0}^2(2u_{2 1}u_{0 2}-u_{0 3}u_{2 0}-u_{1 1}u_{1 2})u_{0 1}+u_{1 0}^3(u_{0 2}u_{1 2}-u_{0 3}u_{1 1})),
\end{split}
\end{equation*}
\begin{equation*}
\begin{split}
I_{0,3}&=\Delta_{2}^{-3}\left(\frac{u_{03}}{3}(u_{01}u_{20}-u_{10}u_{11})^{3}+2(u_{0 1} u_{1 1}-u_{0 2} u_{1 0})(u_{0 1} u_{2 0}-u_{1 0} u_{1 1})\cdot{}\right.\\&\left.\cdot(u_{0 1} u_{1 1} u_{2 1}-u_{0 1} u_{1 2} u_{2 0}-u_{0 2} u_{1 0} u_{2 1}+u_{1 0} u_{1 1} u_{1 2})-\right.{}\\&\left.-\frac{4u_{30}}{3}(u _{0 1}u _{1 1}-u _{0 2}u _{1 0})^3\right).
\end{split}
\end{equation*}

\end{itemize}
\end{example}
\subsection{Weights}
Consider the vector field $V=x\partial_{x}+y\partial_{y}$. Its flow is the scale transformations on the plane $\mathbb{C}^{2}$, and its $\infty$-th prolongation is
\begin{equation*}
V_{\ast}=x\partial_{x}+y\partial_{y}-\sum\limits_{k=1}k\sum\limits_{i=1}^{k}u_{i,k-i}\partial_{u_{i,k-i}}.
\end{equation*}
The vector field $V$, as well as $V_{\ast}$ commutes with the $\mathrm{SL}_{2}(\mathbb{C})$-action and therefore for every $\mathrm{SL}_{2}$-invariant $I$ the function $V_{\ast}(I)$ is invariant too.

We say that invariant $I$ has \textit{weight} $w(I)\in\mathbb{Z}$, if
\begin{equation*}
L_{V_{\ast}}(I)=w(I)I,
\end{equation*}
where $L_{V_{\ast}}$ is the Lie derivative along the vector field $V_{\ast}$.
\begin{example}
\begin{equation*}
w(u_{ij})=-(i+j),\quad w(x)=1,\quad w(\Delta_{2})=-4.
\end{equation*}
\end{example}
Since tensors $\Theta_{k}$ are invariants of affine transformations, $w(\Theta_{k})=0$. Moreover, $w(\omega_{1})=2$, $w(\omega_{2})=0$, and therefore $w(I_{i,j})=-2i$.

Weights can be used to find rational $\mathrm{GL}_{2}(\mathbb{C})$-invariants from polynomial $\mathrm{SL}_{2}(\mathbb{C})$-invariants using the following observation.
\begin{lemma}
Rational $\mathrm{GL}_{2}(\mathbb{C})$-invariants (algebraic or differential) have the form
\begin{equation*}
I=\frac{P}{Q},
\end{equation*}
where $P$ and $Q$ are polynomial $\mathrm{SL}_{2}(\mathbb{C})$-invariants (algebraic or differential) of the same weight.
\end{lemma}
We leave the proof of this lemma to the reader as an exercise.
\subsection{Invariants of binary forms for $n=2,3,4$}
Recall that $\mathcal{B}_{n}\simeq\mathbb{C}^{n+1}$, and the dimension of the group $\mathrm{SL}_{2}(\mathbb{C})$ equals $3$, therefore general orbits have dimension 3 and codimension $(n-2)$, when $n\ge 3$.

An orbit $\mathrm{SL}_{2}(\mathbb{C})\phi$ is said to be \textit{regular}, if the corresponding point on the quotient $\mathbb{C}^{n+1}/\mathrm{SL}_{2}(\mathbb{C})$ is smooth, i.e. there exist $(n-2)$ independent (in a neighborhood of the point) rational invariants $I_{1},\ldots,I_{n-2}$, such that the orbit is given by equations $I_{1}=c_{1},\ldots,I_{n-2}=c_{n-2}$, where $c_{i}$ are constants. Independence means that $dI_{1}\wedge\ldots\wedge dI_{n-2}\ne0$ in the neighborhood of the orbit. Thus $I_{1},\ldots,I_{n-2}$ are regarded as local coordinates on the quotient, and $c_{1},\ldots,c_{n-2}$ are coordinates of the orbit. The Rosenlicht theorem states that all other rational invariants are rational functions of $I_{1},\ldots,I_{n-2}$.

For quadrics $(n=2)$ we have only one differential invariant $\Delta_{2}=u_{20}u_{02}-u_{11}^{2}$. Recall that by replacing $u_{ij}$ with $b_{ij}$ we get algebraic invariants.

For cubics $(n=3)$ we need only $\dim\left(\mathbb{C}^{4}/\mathrm{SL}_{2}(\mathbb{C})\right)=1$ algebraic invariant, which is the discriminant $\Delta_{3}$ of the cubic, and $\dim\left(\mathcal{E}_{3}/\mathrm{SL}_{2}(\mathbb{C})\right)=3$ independent rational differential invariants, which are
\begin{equation}
\label{inv-cub}
J_{1}=\Delta_{2}=u_{02}u_{20}-u_{11}^{2},\quad J_{2}=\nabla_{1}(\Delta_{2}),\quad J_{3}=\Delta_{2}\nabla_{2}(u_{00}).
\end{equation}
Let us restrict differential invariants \eqref{inv-cub} to the cubic $\phi$. We get three functions $J_{1}^{\phi}, J_{2}^{\phi}, J_{3}^{\phi}$ on a plane, namely, binary forms of degrees 2,3,4, therefore, there is one polynomial relation between them:
\begin{equation}
\label{quo-3}
(J_{1}^{\phi})^{5}+(J_{2}^{\phi})^{2}(J_{1}^{\phi})^{2}-16\Delta_{3}(\phi)(J_{3}^{\phi})^{2}=0,
\end{equation}
where $\Delta_{3}(\phi)=\mathrm{Discr}(\phi)$ is the discriminant of the cubic.

Syzygy \eqref{quo-3} can be obtained in Maple using the following code:
\begin{verbatim}
restart;
with(DifferentialGeometry):with(Groebner):
DifferentialGeometry:-Preferences("JetNotation", "JetNotation2"):
with( JetCalculus ):
DGsetup( [x, y], [u], M, 4):
Delta2:=u[0,2]*u[2,0]-u[1,1]^2:
Define invariant derivations according to (7)-(8)
nabla1:=f->u[0,1]*TotalDiff(f,x)-u[1,0]*TotalDiff(f,y):
nabla2:=f->2*(u[0,2]*u[1,0]-u[1,1]*u[0,1])/Delta2*TotalDiff(f,x)+
           2*(u[2,0]*u[0,1]-u[1,1]*u[1,0])/Delta2*TotalDiff(f,y):
Let phi be a binary 3-form
phi:=add(b[i,3-i]*x^i/(i!)*y^(3-i)/(3-i)!,i=0..3):
First invariant (Hessian)
J1:=u[0,2]*u[2,0]-u[1,1]^2:
Second invariant
J2:=nabla1(J1):
Third invariant
J3:=simplify(Delta2*nabla2(u[0,0])):
Restricting invariants to the cubic
Restr:=(f1,f2)->eval(f1,{u[0,0]=f2,
u[0,1]=diff(f2,y),
u[1,0]=diff(f2,x),
u[2,0]=diff(f2,x$2),
u[0,2]=diff(f2,y$2),
u[1,1]=diff(f2,[x,y]),
u[3,0]=diff(f2,x$3),
u[2,1]=diff(f2,[x,x,y]),
u[1,2]=diff(f2,[x,y,y]),
u[0,3]=diff(f2,y$3)}):
Restriction of J1 to the cubic
J1phi:=Restr(J1,phi):
Restriction of J2 to the cubic
J2phi:=Restr(J2,phi):
Restriction of J3 to the cubic
J3phi:=Restr(J3,phi):
Finding syzygy
syz1:=Basis([J1phi-Z0, J2phi-Z2, J3phi-Z3],plex(x, y, Z0, Z2, Z3))[1]:
\end{verbatim}
Removing the restriction to the cubic $\phi$ from \eqref{quo-3}, we get a differential equation of the third order:
\begin{equation}
\label{quo-3-1}
\left\{(J_{1})^{5}+(J_{2})^{2}(J_{1})^{2}-16\Delta_{3}(\phi)(J_{3})^{2}=0\right\}\subset\mathbf{J}^{3}.
\end{equation}
Thus we have the following criterion of $\mathrm{SL}_{2}(\mathbb{C})$-equivalence of binary 3-forms:
\begin{theorem}
Let $\phi$ be a regular binary 3-form ($\Delta_{3}(\phi)\ne0$). Then, $\mathrm{SL}_{2}(\mathbb{C})$-orbit of $\phi$ consists of solutions to the third order differential equation \eqref{quo-3-1} together with $\mathcal{E}_{3}$.
\end{theorem}

For quartics $(n=4)$ we take the following differential invariants
\begin{equation*}
J_{0}=u_{00},\quad J_{2}=\Delta_{2}=u_{02}u_{20}-u_{11}^{2},\quad J_{3}=-\nabla_{1}(J_{2}).
\end{equation*}
Again, if we restrict these invariants to a regular quartic $\phi$, we will obtain quartics $J_{0}^{\phi}$, $J_{2}^{\phi}$, $J_{3}^{\phi}$ on the plane, and the polynomial relation between them is
\begin{equation}
\label{quo-4}
9(J_{3}^{\phi})^{2}+16(J_{2}^{\phi})^{3}+144\alpha(J_{0}^{\phi})^{2}J_{2}^{\phi}+864\delta(J_{0}^{\phi})^{3}=0,
\end{equation}
where
\begin{equation*}
\alpha=4b_{13}b_{31}-b_{40}b_{04}-3b_{22}^{2}
\end{equation*}
is the Hankel apolar, and
\begin{equation*}
\delta=b_{22}b_{40}b_{04}-b_{04}b_{31}^{2}-b_{40}b_{13}^{2}+2b_{13}b_{22}b_{31}-b_{22}^{3}
\end{equation*}
is the Hankel determinant.

Relation \eqref{quo-4} can be obtained by means of the same Maple code as we used for cubics.

Removing the restriction to the quartic $\phi$ from \eqref{quo-4}, we get a differential equation of the third order:
\begin{equation}
\label{quo-4-1}
\left\{9(J_{3})^{2}+16(J_{2})^{3}+144\alpha(J_{0})^{2}J_{2}+864\delta(J_{0})^{3}=0\right\}\subset\mathbf{J}^{3}.
\end{equation}

Thus we have a similar theorem for quartics:
\begin{theorem}
Let $\phi$ be a regular binary 4-form. Then, $\mathrm{SL}_{2}(\mathbb{C})$-orbit of $\phi$ consists of solutions to the third order differential equation \eqref{quo-4-1} together with $\mathcal{E}_{4}$.
\end{theorem}
\section{Quotients}
\label{sec:3}
This section gives a general introduction into the structure of quotients of algebraic manifolds and equations under the action of algebraic groups. The main results are given by the Rosenlicht and the Lie-Tresse theorems.
\subsection{Rosenlicht theorem}
Let $\Omega$ be a set with an action of a group $G$:
\begin{equation*}
G\times\Omega\to\Omega,\quad g\times\omega\mapsto g\omega,
\end{equation*}
Then, the set $G/\Omega$ of all $G$-orbits is called \textit{quotient}:
\begin{equation*}
\Omega/G=\bigcup\limits_{\omega\in\Omega}\left\{G\omega\right\}.
\end{equation*}
\begin{remark}
The projection $\pi\colon\Omega\to\Omega/G$ allows us to identify functions on the quotient $\Omega/G$ with functions on $\Omega$ that are $G$-invariants, i.e. $f\circ g=f$.
\end{remark}
Let $\Omega$ be a topological space, $G$ be a topological group and let $G$-action be continuous. Then, the quotient $\Omega/G$ is naturally a topological space, that is, a subset $U\subset\Omega/G$ is said to be open if and only if the preimage $\pi^{-1}(U)\subset\Omega$ is open.
\begin{remark}
In general, we cannot guarantee that the quotient $\Omega/G$ shall inherit topological properties (e.g. the Hausdorff condition) of $\Omega$.
\end{remark}
\begin{example}
\begin{enumerate}
\item Let $\Omega=\mathbb{R}^{2}$, $G=\mathrm{SL}_{2}(\mathbb{R})$, and $\mathrm{SL}_{2}(\mathbb{R})\times\mathbb{R}^{2}\to\mathbb{R}^{2}$ be the natural action. Then,
\begin{equation*}
\mathbb{R}^{2}/\mathrm{SL}_{2}(\mathbb{R})=\mathbf{0}\cup\bigstar,
\end{equation*}
where $\mathbf{0}=\mathrm{SL}_{2}(\mathbb{R})(0)$ is the orbit of the origin, $0\in\mathbb{R}^{2}$, and $\bigstar$ is the orbit of any nonzero point. This is an example of the famous \textit{Sierpinski topological space}, consisting of two points, one of which $\mathbf{0}$ is closed, but another one $\bigstar$ is open.
\item Let $\Omega=\mathbb{R}^{2}$, $G=\mathbb{R}^{\ast}=\mathbb{R}\setminus0$, and $\mathbb{R}^{\ast}\times\mathbb{R}^{2}\to\mathbb{R}^{2}$ be the natural action. Then,
\begin{equation*}
\mathbb{R}^{2}/\mathbb{R}^{\ast}=\mathbf{0}\cup\mathbb{R}P^{1},
\end{equation*}
where $\mathbb{R}P^{1}$ is the projective 1-dimensional space.
\end{enumerate}
\end{example}
If $\Omega$ is a smooth manifold and $G$ is a Lie group, then we have no way to determine whether the quotient $\Omega/G$ is also a smooth manifold, except for the case when $G$-action is free and proper.

Let $G$ be an algebraic manifold (an irreducible variety without singularities over a field of zero characteristic), $G$ be an algebraic group, and $G\times\Omega\to\Omega$ be an algebraic action. By $\mathcal{F}(\Omega)$ we denote the field of rational functions on $\Omega$ and by $\mathcal{F}(\Omega)^{G}\subset\mathcal{F}(\Omega)$ the field of rational $G$-invariants. An orbit $G\omega\subset\Omega$ (as well as the point $\omega$) is said to be \textit{regular}, if there are $m=\mathrm{codim}(G\omega)$ $G$-invariants $x_{1},\ldots,x_{m}$, such that their differentials are linear independent at the points of the orbit.

Let $\Omega_{0}=\Omega\setminus\mathrm{Sing}$ be the set of all regular points and $Q(\Omega)=\Omega_{0}/G$ be the set of all regular orbits.
\begin{theorem}[Rosenlicht, \cite{VinPop,Ros}]
The set $\Omega_{0}$ is open and dense in $\Omega$ in the Zariski topology.
\end{theorem}
Invariants $x_{1},\ldots,x_{m}$ can be considered as local coordinates on the quotient $Q(\Omega)$ in the neighborhood of the point $G\omega\in Q(\Omega)$. On intersections of charts these coordinates are related by rational functions, which means that $Q(\Omega)$ is an algebraic manifold of the dimension $m=\mathrm{codim}(G\omega)$. Thus we have the rational map $\pi\colon\Omega_{0}\to Q(\Omega)$ of algebraic manifolds, which gives us a field isomorphism $\mathcal{F}(\Omega)^{G}=\pi^{\ast}(\mathcal{F}(Q(\Omega)))$.

It is essential that the Rosenlicht's theorem is valid only for algebraic manifolds. Indeed, following the algebraic case, let $\Omega$ be a smooth manifold, and $G$ be a Lie group. An orbit $G\omega$ (as the point $\omega$ itself) is said to be \textit{regular}, if there are $m=\mathrm{codim}(G\omega)$ smooth independent (in the above sense) invariants. Again, let $\Omega_{\mathrm{reg}}\subset\Omega$ be the set of regular points, then the quotient $\Omega_{\mathrm{reg}}/G$ is a smooth manifold, and the projection $\pi\colon\Omega_{\mathrm{reg}}\to\Omega_{\mathrm{reg}}/G$ gives us an isomorphism of algebras $C^{\infty}(\Omega_{\mathrm{reg}})^{G}$ and $C^{\infty}(\Omega_{\mathrm{reg}}/G)$, $\pi^{\ast}\left(C^{\infty}(\Omega_{\mathrm{reg}}/G)\right)=C^{\infty}(\Omega_{\mathrm{reg}})^{G}$. In contrast to the algebraic case we could not guarantee that $\Omega_{\mathrm{reg}}$ is dense in $\Omega$.

Let, again, $\Omega$ be an algebraic manifold, and let $\mathfrak{g}$ be a Lie subalgebra of the Lie algebra of vector fields on $\Omega$. The Lie algebra $\mathfrak{g}$ is said to be \textit{algebraic} if there exists an algebraic action of the algebraic group $G$, such that $\mathfrak{g}$ coincides with the image of the Lie algebra $\mathrm{Lie}(G)$ under this action. By an \textit{algebraic closure} of the Lie algebra $\mathfrak{g}$ we mean an intersection of all algebraic Lie algebras, containing $\mathfrak{g}$.
\begin{example}
\begin{enumerate}
\item $\Omega=\mathbb{R}$, the Lie algebra
\begin{equation*}
\mathfrak{g}=\mathfrak{sl}_{2}=\langle\partial_{x},x\partial_{x},x^{2}\partial_{x}\rangle
\end{equation*}
is algebraic.
\item $\Omega=\mathbb{R}^{2}$, and the Lie algebra
\begin{equation*}
\mathfrak{g}=\langle x\partial_{x}+\lambda y\partial_{y}\rangle
\end{equation*}
is algebraic if $\lambda\in\mathbb{Q}$. In the case $\lambda\notin\mathbb{Q}$ the closure is $\tilde{\mathfrak{g}}=\langle x\partial_{x},y\partial_{y}\rangle$.
\item $\Omega=S^{1}\times S^{1}$ --- torus, the Lie algebra
\begin{equation*}
\mathfrak{g}=\langle\partial_{\phi}+\lambda\partial_{\psi}\rangle
\end{equation*}
is algebraic if $\lambda\in\mathbb{Q}$. In the case $\lambda\notin\mathbb{Q}$ the closure is $\tilde{\mathfrak{g}}=\langle \partial_{\phi},\partial_{\psi}\rangle$.
\end{enumerate}
\end{example}

It turns out that the Rosenlicht theorem is also valid for algebraic Lie algebras, or for algebraic closure in the case of general Lie algebras.

Indeed, let $\mathfrak{g}$ be a Lie algebra of vector fields on an algebraic manifold $\Omega$ and let $\tilde{\mathfrak{g}}$ be its algebraic closure. Then, the field $\mathcal{F}(\Omega)^{\mathfrak{g}}$ of rational $\mathfrak{g}$-invariants has a transcendence degree equal to the codimension of  $\tilde{\mathfrak{g}}$-orbits that is the dimension of the quotient $Q(\Omega)$.
\subsection{Algebraicity in Jet Geometry}
Let $\pi\colon E(\pi)\to M$ be a smooth bundle over a manifold $M$ and let $\pi_{k}\colon\mathbf{J}^{k}\to M$ be the bundle of sections of $k$-jets.

The manifold $\mathbf{J}^{k}$ is equipped with the Cartan distribution, which in canonical jet coordinates $(x,u^{j}_{\sigma})$ is given by differential 1-forms
\begin{equation}
\label{CartD}
\varkappa^{j}_{\sigma}=du^{j}_{\sigma}-\sum\limits_{i}u^{j}_{\sigma i}dx_{i}.
\end{equation}
The Lie-B\"{a}klund theorem \cite{KrVin,KrVinLych} states that types of Lie transformations, i.e. local diffeomorphisms of $\mathbf{J}^{k}$ preserving the Cartan distribution \eqref{CartD}, are determined by the dimension of $\pi$, namely, they are prolongations of
\begin{itemize}
\item the pseudogroup $\mathrm{Cont}(\pi)$ of local \textit{contact transformations} of $\mathbf{J}^{1}$, in the case $\dim\pi=1$;
\item the pseudogroup $\mathrm{Point}(\pi)$ of local \textit{point transformations} of $\mathbf{J}^{0}$, i.e. local diffeomorphisms of $\mathbf{J}^{0}$, in the case $\dim\pi>1$.
\end{itemize}
Moreover, it is known that
\begin{itemize}
\item all bundles $\pi_{k,k-1}\colon\mathbf{J}^{k}\to\mathbf{J}^{k-1}$ are affine bundles for $k\ge2$, when $\dim\pi\ge2$, and for $k\ge3$, when $\dim\pi=1$;
\item prolongations of pseudogroups in canonical jet coordinates $(x,u^{j}_{\sigma})$ are given by rational in $u^{j}_{\sigma}$ functions.
\end{itemize}
Therefore,
\begin{itemize}
\item in the case $\dim\pi\ge2$ the fibres $\mathbf{J}_{\theta}^{k,0}$ of the projections $\pi_{k,0}\colon\mathbf{J}^{k}\to\mathbf{J}^{0}$ at points $\theta\in\mathbf{J}^{0}$ are algebraic manifolds, and the stationary subgroup $\mathrm{Point}_{\theta}(\pi)\subset\mathrm{Point}(\pi)$ gives us birational isomorphisms of the manifold;
\item in the case $\dim\pi=1$ the fibres $\mathbf{J}_{\theta}^{k,1}$ of the projections $\pi_{k,1}\colon\mathbf{J}^{k}\to\mathbf{J}^{1}$ at points $\theta\in\mathbf{J}^{1}$ are algebraic manifolds, and the stationary subgroup $\mathrm{Cont}_{\theta}(\pi)\subset\mathrm{Cont}(\pi)$ gives us birational isomorphisms of the manifold.
\end{itemize}
\subsection{Algebraic Differential Equations}
A differential equation $\mathcal{E}_{k}\subset\mathbf{J}^{k}$ is said to be \textit{algebraic}, if fibres $\mathcal{E}_{k,\theta}$ of the projections $\pi_{k,0}\colon\mathcal{E}_{k}\to\mathbf{J}^{0}$, when $\dim\pi\ge2$, or $\pi_{k,1}\colon\mathcal{E}_{k}\to\mathbf{J}^{1}$, when $\dim\pi=1$ , are algebraic manifolds.
\begin{remark}
If $\mathcal{E}_{k}$ is algebraic and formally integrable, then the prolongations $\mathcal{E}_{k}^{(l)}=\mathcal{E}_{k+l}\subset\mathbf{J}^{k+l}$ are algebraic too.
\end{remark}
By a symmetry algebra of algebraic differential equations we mean one of the following:
\begin{itemize}
\item for $\dim\pi\ge2$, a Lie algebra $\mathrm{sym}(\mathcal{E}_{k})$ of point symmetries (point vector fields), which is transitive on $\mathbf{J}^{0}$, and stationary subalgebras $\mathrm{sym}_{\theta}(\mathcal{E}_{k})$, $\theta\in\mathbf{J}^{0}$, produce actions of algebraic Lie algebras on algebraic manifolds $\mathcal{E}_{l,\theta}$, for all $l\ge k$;
\item for $\dim\pi=1$, a Lie algebra $\mathrm{sym}(\mathcal{E}_{k})$ of contact symmetries (contact vector fields), which is transitive on $\mathbf{J}^{1}$, and stationary subalgebras $\mathrm{sym}_{\theta}(\mathcal{E}_{k})$, $\theta\in\mathbf{J}^{1}$, produce actions of algebraic Lie algebras on algebraic manifolds $\mathcal{E}_{l,\theta}$, for all $l\ge k$.
\end{itemize}

Let $\mathcal{E}_{k}$ be a formally integrable algebraic differential equation,  $\mathcal{E}_{l}$ be its $(l-k)$-prolongation, and $\mathfrak{g}$ be its algebraic symmetry Lie algebra. Then, all the $\mathcal{E}_{l}$ are algebraic manifolds, and we have a tower of algebraic bundles:
\begin{equation*}
\mathcal{E}_{k}\longleftarrow\mathcal{E}_{k+1}\longleftarrow\cdots\longleftarrow \mathcal{E}_{l}\longleftarrow \mathcal{E}_{l+1}\longleftarrow\cdots.
  \end{equation*}
  A point $\theta\in\mathcal{E}_{l}$ (a $\mathfrak{g}$-orbit) is said to be \textit{strongly regular}, if it is regular and its projection to $\mathcal{E}_{l-i}$ for all $i=1,...,l-k$ is regular too.

Let $\mathcal{E}_{l}^{0}\subset\mathcal{E}_{l}$ be the set of all strongly regular points and $Q_{l}(\mathcal{E})$ be the set of all regular $\mathfrak{g}$-orbits. Then, due to the Rosenlicht's theorem,  $Q_{l}(\mathcal{E})$ are algebraic manifolds, and projections $\varkappa_{l}\colon\mathcal{E}_{l}^{0}\to Q_{l}(\mathcal{E})$ are rational maps, such that $\varkappa_{l}^{\ast}(\mathcal{F}(Q_{l}(\mathcal{E})))=\mathcal{F}(\mathcal{E}_{l}^{0})^{\mathfrak{g}}$, where $\mathcal{F}(Q_{l}(\mathcal{E}))$ is the field of rational functions on $Q_{l}(\mathcal{E})$, and $\mathcal{F}(\mathcal{E}_{l}^{0})^{\mathfrak{g}}$ is the field of rational $\mathfrak{g}$-invariant functions (\textit{rational differential invariants}).

Since the $\mathfrak{g}$-action preserves the Cartan distribution $\mathcal{C}(\mathcal{E}_{l})$, projections $\varkappa_{l}$ define distributions on the quotients $Q_{l}(\mathcal{E})$. Finally, we have the tower of algebraic bundles of the quotients
\begin{equation}
\label{seq-quo}
Q_{k}(\mathcal{E})\stackrel{\pi_{k+1,k}}{\longleftarrow} Q_{k+1}(\mathcal{E})\longleftarrow\cdots\longleftarrow Q_{l}(\mathcal{E})\stackrel{\pi_{l+1,l}}{\longleftarrow}Q_{l+1}(\mathcal{E})\longleftarrow\cdots,
  \end{equation}
such that $(\pi_{l+1,l})_{\ast}(\mathcal{C}(Q_{l+1}(\mathcal{E})))=\mathcal{C}(Q_{l}(\mathcal{E}))$ for $l\ge k$.

Locally, sequence \eqref{seq-quo} has the same structure as for some equation $F$, which is called a \textit{quotient PDE}.

\subsection{Lie-Tresse theorem}
First, we discuss Lie-Tresse derivatives, which are necessary for description of quotient PDEs.

Let $\omega\in\Omega^{1}(\mathbf{J}^{k})$ be a differential 1-form on the space of $k$-jets and let $\mathcal{C}_{k}$ be the Cartan distribution. Then, the class
\begin{equation*}
\omega^{h}=\pi_{k+1,k}^{\ast}(\omega)\quad\mathrm{mod}\,\mathrm{Ann}(\mathcal{C}_{k+1})
\end{equation*}
is called a \textit{horizontal part} of $\omega$. In the canonical jet coordinates $(x,u^{j}_{\sigma})$ we have
\begin{equation*}
\omega=\sum\limits_{i=1}^{n}a_{i}dx_{i}+\sum\limits_{\substack{j\le m\\|\sigma|\le k}}b^{j}_{\sigma}du^{j}_{\sigma},
\end{equation*}
and its horizontal part is
\begin{equation*}
\omega^{h}=\sum\limits_{\substack{j\le m\\|\sigma|\le k\\i\le n}}\left(a_{i}+b^{j}_{\sigma}u^{j}_{\sigma i}\right)dx_{i},
\end{equation*}
where $n=\dim M$, $m=\dim\pi$.

 Applying this construction to the differential $df$ of the function $f\in C^{\infty}(\mathbf{J}^{k})$ we get a \textit{total differential} $\widehat{d}f=(df)^{h}$. In canonical coordinates it is
\begin{equation*}
\widehat{d}f=\sum\limits_{i=1}^{n}\frac{df}{dx_{i}}dx_{i},\quad\frac{d}{dx_{i}}=\frac{\partial}{\partial x_{i}}+\sum\limits_{j,\sigma}u^{j}_{\sigma i}\frac{\partial}{\partial u^{j}_{\sigma}}.
\end{equation*}
It is worth mentioning that the operation of taking the horizontal part as well as total differentials are invariant with respect to point and contact transformations.

Functions $f_{1},\ldots,f_{n}\in C^{\infty}(\mathbb{J}^{k})$ are said to be \textit{in general position} in some domain $D$ if
\begin{equation}
\label{indep}
\widehat{d}f_{1}\wedge\ldots\wedge\widehat{d}f_{n}\ne0\text{ in }D.
\end{equation}
Given fixed $f_{1},\ldots,f_{n}$ satisfying \eqref{indep} one has the following decomposition for $f\in C^{\infty}(\mathbf{J}^{k})$ in $D$:
\begin{equation*}
\widehat{d}f=\sum\limits_{i=1}^{n}F_{i}\widehat{d}f_{i},
\end{equation*}
where $F_{i}$ are smooth functions in the domain $\pi^{-1}_{k+1,k}(D)\subset\mathbf{J}^{k+1}$, called \textit{Tresse derivatives} and denoted by $\frac{df}{df_{i}}$.
\begin{theorem}
Let $f_{1},\ldots,f_{n}$ be $\mathfrak{g}$-invariants of order $\le k$ in general position. Then, for any $\mathfrak{g}$-invariant $f$ of order $\le k$ the Tresse derivatives $\frac{df}{df_{i}}$ are $\mathfrak{g}$-invariants of order $\le k+1$.
\end{theorem}
\begin{example}
Consider the action of the Lie group of translations on a plane. Its Lie algebra is
\begin{equation*}
\mathfrak{g}=\langle\partial_{x},\partial_{y}\rangle.
\end{equation*}
Let us take its invariants $f_{1}=u_{00}$, $f_{2}=u_{10}$, $f=u_{01}$. Then, the Tresse derivatives are of the form
\begin{eqnarray*}
 \frac{d}{df_{1}} &=&\frac{u_{11}}{u_{10}u_{11}-u_{01}u_{20}}\frac{d}{dx}+\frac{u_{20}}{u_{01}u_{20}-u_{10}u_{11}}\frac{d}{dy}, \\
 \frac{d}{df_{2}} &=&\frac{u_{01}}{u_{01}u_{20}-u_{10}u_{11}}\frac{d}{dx}+\frac{u_{10}}{u_{10}u_{11}-u_{01}u_{20}}\frac{d}{dy}.
\end{eqnarray*}
Applying them to the differential invariant $f=u_{01}$ of the first order, we get two more invariants of the second order:
\begin{equation*}
J_{1}=\frac{df}{df_{1}}=\frac{u_{20}u_{02}-u_{11}^{2}}{u_{10}u_{20}-u_{10}u_{11}},\quad J_{2}=\frac{df}{df_{2}}=\frac{u_{01}u_{11}-u_{02}u_{10}}{u_{01}u_{20}-u_{10}u_{11}}.
\end{equation*}
\end{example}
The following statement known as the \textit{global Lie-Tresse theorem} \cite{KrugLych} gives the conditions of finiteness for a generating set of invariants of a pseudogroup action on a differential equation:
\begin{theorem}[Kruglikov, Lychagin]
Let $\mathcal{E}_{k}\subset\mathbf{J}^{k}$ be an algebraic formally integrable differential equation and let $\mathfrak{g}$ be its algebraic symmetry Lie algebra. Then, there exist rational differential $\mathfrak{g}$-invariants $a_{1},\ldots,a_{n},b^{1},\ldots,b^{N}$ of order $\le l$, such that the field of rational $\mathfrak{g}$-invariants is generated by rational functions of these functions and Tresse derivatives $\frac{d^{|\alpha|}b^{j}}{da^{\alpha}}$.
\end{theorem}
Local version of this result goes back to S. Lie and A. Tresse.
\begin{remark}
\begin{enumerate}
\item In contrast to algebraic invariants, where we have only algebraic operations, in the case of differential invariants we have more operations. Namely, the Tresse derivatives give us new differential invariants.
\item The algebra of differential invariants is not freely generated, there are relations between invariants, called \textit{syzygies}. The syzygies provide us with new differential equations, called \textit{quotient equations}.
\item From the geometrical viewpoint, the Lie-Tresse theorem states that there is a level $l$ and a domain $D\subset Q_{l}(\mathcal{E})$, where invariants $a_{1},\ldots,a_{n},b^{1},\ldots,b^{N}$ serve as local coordinates, and the preimage of $D$ in the tower
\begin{equation}
\label{tower}
Q_{l}(\mathcal{E})\stackrel{\pi_{l+1,l}}{\longleftarrow} Q_{l+1}(\mathcal{E})\longleftarrow\cdots\longleftarrow Q_{r}(\mathcal{E})\stackrel{\pi_{r+1,r}}{\longleftarrow}Q_{r+1}(\mathcal{E})\longleftarrow\cdots
  \end{equation}
  is an infinitely prolonged differential equation given by the syzygy. For this reason we call the quotient tower \eqref{tower} an \textit{algebraic diffiety}.
\end{enumerate}
\end{remark}
\subsection{Integrability via Quotients}
Here we discuss the importance of above constructions for integrability of differential equations. First, let us summarize the relations between differential equations and their quotients:
\begin{enumerate}
\item Let $L$ be a solution to a differential equation $\mathcal{E}$ (in the sense of integral manifolds of the Cartan distribution) and let $a_{i}|_{L}, b^{j}|_{L}$ be the values of differential invariants on the solution $L$. Then, we have $b^{j}|_{L}=B^{j}\left(a|_{L}\right)$, and functions $B^{j}$ are exactly solutions to the quotient differential equations.
\item Let $b^{j}=B^{j}(a)$ be a solution to a quotient PDE. Then, adding differential constraints $b^{j}-B^{j}(a)=0$ we get a finite type equation $\mathcal{E}\cap\left\{b^{j}-B^{j}(a)=0\right\}$ with solutions being a $\mathfrak{g}$-orbit of a solution to $\mathcal{E}$. This gives us a method of finding compatible constraints to be added to the original system of PDEs, which reduces the integration of the PDE to the integration of a completely integrable Cartan distribution having the same symmetry algebra. This is essential for finding smooth solutions, as well as those with singularities \cite{Eiv,LR-lob}.
\item Symmetries of quotient PDEs are B\"{a}cklund-type transformations for the equation $\mathcal{E}$.
\end{enumerate}
Let us now illustrate this on examples. As an exercise, we recommend the reader to do the computations for these examples.
\begin{example}
\begin{enumerate}
\item Invariants of the Lie algebra $\mathfrak{g}=\langle\partial_{x}\rangle$ of $x$-translations on the line $\Omega=\mathbb{R}$ are generated by
\begin{equation*}
\langle a=u_{0}, b=u_{1}\rangle
\end{equation*}
and Tresse derivative
\begin{equation*}
\frac{d}{da}=u_{1}^{-1}\frac{d}{dx}.
\end{equation*}
Then, for the $x$-invariant ODE of the third order $F(u_{0},u_{1},u_{2},u_{3})=0$ the quotient equation is of order 2 and has the form
\begin{equation*}
F\left(a,b,b\frac{db}{da},b^{2}\frac{d^{2}b}{da^{2}}\right)=0.
\end{equation*}
This is a standard reduction of order for ODEs of the form $F(u_{0},u_{1},u_{2},u_{3})=0$.

Let us now choose other Lie-Tresse coordinates: \begin{equation*}\langle a=u_{2},b^{1}=u_{0},b^{2}=u_{1}\rangle\end{equation*} and Tresse derivative
\begin{equation*}
\frac{d}{da}=u_{3}^{-1}\frac{d}{dx}.
\end{equation*}
In this case, the quotient equation for $F(u_{0},u_{1},u_{2},u_{3})=0$ is a system of ODEs:
\begin{equation*}
F\left(b^{1},b^{2},a,a\left(\frac{db^{2}}{da}\right)^{-1}\right)=0,\quad a\frac{db^{1}}{da}-b^{2}\frac{db^{2}}{da}=0.
\end{equation*}
\item Invariants of the Lie algebra $\mathfrak{g}=\langle\partial_{x},x\partial_{x}\rangle$ of affine transformations of the line $\Omega=\mathbb{R}$ are \begin{equation*}
\left\langle u_{0},\frac{u_{2}}{u_{1}^{2}},\frac{u_{3}}{u_{1}^{3}},\frac{u_{4}}{u_{1}^{4}},\ldots\right\rangle.
\end{equation*}
Let us take \begin{equation*}\left\langle a=u_{0},b=\frac{u_{2}}{u_{1}^{2}}\right\rangle\end{equation*}
and consider a $\mathfrak{g}$-invariant equation
\begin{equation*}
F\left(u_{0},\frac{u_{2}}{u_{1}^{2}},\frac{u_{3}}{u_{1}^{3}},\frac{u_{4}}{u_{1}^{4}}\right)=0.
\end{equation*}
Its quotient will be
\begin{equation*}
F\left(a,b,\frac{db}{da}+2b^{2},\frac{d^{2}b}{da^{2}}+6b\frac{db}{da}+6b^{3}\right)=0.
\end{equation*}
\item Invariants of the Lie algebra $\mathfrak{g}=\mathfrak{sl}_{2}(\mathbb{R})=\langle\partial_{x},x\partial_{x},x^{2}\partial_{x}\rangle$ on the line $\Omega=\mathbb{R}$ are
\begin{equation*}
\left\langle u_{0},\frac{u_{3}}{u_{1}^{3}}-\frac{3u_{2}^{2}}{2u_{1}^{4}},\frac{u_{4}}{u_{1}^{4}}-6\frac{u_{2}u_{3}}{u_{1}^{5}}+6\frac{u_{2}^{3}}{u_{1}^{6}},\ldots\right\rangle.
\end{equation*}
Let us take \begin{equation*}\left\langle a=u_{0},b=\frac{u_{3}}{u_{1}^{3}}-\frac{3u_{2}^{2}}{2u_{1}^{4}}\right\rangle\end{equation*}
and consider a $\mathfrak{g}$-invariant equation
\begin{equation*}
F\left( u_{0},\frac{u_{3}}{u_{1}^{3}}-\frac{3u_{2}^{2}}{2u_{1}^{4}},\frac{u_{4}}{u_{1}^{4}}-6\frac{u_{2}u_{3}}{u_{1}^{5}}+6\frac{u_{2}^{3}}{u_{1}^{6}}\right)=0.
\end{equation*}
Its quotient will be
\begin{equation*}
F\left(a,b,\frac{db}{da}\right)=0.
\end{equation*}
\item Invariants of the Lie algebra $\mathfrak{g}=\langle\partial_{x},\partial_{y}\rangle$ on the plane $\Omega=\mathbb{R}^{2}$ are
\begin{equation*}
\left\langle u_{00},u_{10},u_{01},u_{20},u_{11},u_{02}\ldots\right\rangle.
\end{equation*}
Let us take \begin{equation*}\left\langle a_{1}=u_{10},a_{2}=u_{01},b^{1}=u_{00},b^{2}=u_{11}\right\rangle\end{equation*} as Lie-Tresse coordinates. Then, assuming $b^{1}=B^{1}(a_{1},a_{2})$, $b^{2}=B^{2}(a_{1},a_{2})$, we have
\begin{equation*}
B^{1}_{a_{1}}=\delta^{-1}(u_{10}u_{02}-u_{01}u_{11}),\quad B^{1}_{a_{2}}=\delta^{-1}(u_{01}u_{20}-u_{10}u_{11}),
\end{equation*}
\begin{equation*}
B^{2}_{a_{1}}=\delta^{-1}(u_{02}u_{21}-u_{11}u_{12}),\quad B^{2}_{a_{2}}=\delta^{-1}(u_{20}u_{12}-u_{11}u_{21}),
\end{equation*}
where $\delta=u_{20}u_{02}-u_{11}^{2}$ is the Hessian. The syzygies
\begin{eqnarray*}
0&=&-B^{1}_{a_{2}a_{2}}B^{2}B^{1}_{a_{1}a_{1}}+B^{2}(B^{1}_{a_{1}a_{2}})^{2}-B^{1}_{a_{1}a_{2}},\\
0&=&a_{1}B^{1}_{a_{1}a_{1}}+a_{2}B^{1}_{a_{1}a_{2}}-B^{1}_{a_{1}},\\
0&=&a_{1}B^{2}B^{1}_{a_{1}a_{1}}B^{1}_{a_{1}a_{2}}+a_{2}B^{2}(B^{1}_{a_{1}a_{2}})^{2}-B^{2}B^{1}_{a_{1}a_{1}}B^{1}_{a_{2}}-a_{2}B^{1}_{a_{1}a_{2}}
\end{eqnarray*}
are quotient PDEs for the equation $u_{11}=B^{2}(u_{10},u_{01})$.

In particular, equation $u_{11}=0$ is \textit{self-dual}, it coincides with its quotient.
\end{enumerate}
\end{example}
\begin{remark}
\begin{enumerate}
\item If an ODE of order $k$ admits a solvable symmetry Lie algebra $\mathfrak{g}$, and $\dim\mathfrak{g}=k$, then the integration can be done explicitly using the Lie-Bianchi theorem. If the Lie algebra $\mathfrak{g}$ is not solvable, but still $\dim\mathfrak{g}=k$, then the integration can be done by means of model equations \cite{LychWisla2019}.
\item If $\dim\mathfrak{g}=k-1$, the integration splits into the integration of the first order quotient equation and integration of $(k-1)$ order equation with the same symmetry algebra $\mathfrak{g}$. Continuing, we reduce the integration to the integration to a series of quotients.
\end{enumerate}
\end{remark}
\section{Algebraic Plane Curves}
\label{sec:4}
This section is devoted to finding affine invariants for algebraic plane curves using affine connections.
\subsection{Connections and Affine Structures}
The motivation to study connections goes back to classical mechanics, when one needs to define acceleration. If we consider a vector field $Y$ on a manifold $M$ as the field of velocities, then we should be able to compare tangent vectors at different points of the manifold. Let $x(t)$ be a path on the manifold $M$ and assume that we have linear isomorphisms $\lambda(t)\colon T_{x(t)}M\to T_{x(0)}M$ of tangent spaces. Then, taking images $Y(t)=\lambda(t)\left(Y_{x(t)}\right)\in T_{x(0)}M$ of vectors $Y(t)\in T_{x(t)}M$, we get the velocity of variation of the vector field along the path $x(t)$:
\begin{equation}
\label{cov-der}
\left.\frac{dY(t)}{dt}\right|_{t=0}\in T_{x(0)}M.
\end{equation}
Let $x(t)$ be the trajectory of another vector field $X$ on the manifold $M$. Then, taking derivatives \eqref{cov-der} at points of $M$, we get a vector field $\nabla_{X}Y$ on $M$. Assuming that the map $X\times Y\to\nabla_{X}Y$ is $C^{\infty}(M)$-linear in $X$, we obtain the notion of a \textit{covariant derivative}.

Let $M$ be a smooth manifold and let $\mathcal{D}(M)$ be the module of vector fields on $M$. Then, the \textit{covariant derivative} is a map
\begin{equation*}
\nabla_{X}\colon\mathcal{D}(M)\to\mathcal{D}(M),\quad X\in\mathcal{D}(M),
\end{equation*}
satisfying conditions
\begin{enumerate}
\item $\nabla_{X_{1}+X_{2}}=\nabla_{X_{1}}+\nabla_{X_{2}}$
\item $\nabla_{fX}=f\nabla_{X},\, f\in C^{\infty}(M)$,
\item $\nabla_{X}(Y_{1}+Y_{2})=\nabla_{X}(Y_{1})+\nabla_{X}(Y_{2})$
\item $\nabla_{X}(fY)=X(f)Y+f\nabla_{X}(Y)$,
\end{enumerate}
where $X_{i}, Y_{i}, X, Y\in\mathcal{D}(M)$, $f\in C^{\infty}(M)$. Any affine (linear) connection on a manifold $M$ is defined by its covariant derivative.

Let $\nabla$ and $\tilde{\nabla}$ be two affine connections, then the difference $\Gamma_{X}=\nabla_{X}-\tilde{\nabla}_{X}\colon\mathcal{D}(M)\to\mathcal{D}(M)$ is a linear operator, $\Gamma_{X}\in\mathrm{End}(\mathcal{D}(M))$, i.e. a map $X\mapsto\Gamma_{X}$ is $\mathbb{R}$-linear, and $\Gamma_{X}(fY)=f\Gamma_{X}(Y)$. In other words, $\Gamma\in\mathrm{End}(\mathcal{D}(M))\otimes\Omega^{1}(M)$ is an $\mathrm{End}(\mathcal{D}(M))$-valued differential one-form on $M$, called \textit{connection form}, and finding connection on a manifold is equivalent to finding a connection form.

Let $M=\mathbb{R}^{n}$ with coordinates $(x_{1},\ldots,x_{n})$ be a real vector space. Consider $M$ as an affine space with standard identifications of tangent spaces at different points, we come to the covariant derivatives
\begin{equation*}
\nabla^{s}_{\partial_{i}}(\partial_{j})=0,
\end{equation*}
and any other connection has the form
\begin{equation*}
\nabla_{\partial_{i}}(\partial_{j})=\sum\limits_{k}\Gamma_{ij}^{k}\partial_{k},
\end{equation*}
where now and further on $\partial_{i}=\partial_{x_{i}}$, $d_{i}=dx_{i}$, $\Gamma_{ij}^{k}$ are Christoffel symbols.

The \textit{torsion tensor} $T$ of a connection $\nabla$ is
\begin{equation*}
T(X,Y)=\nabla_{X}(Y)-\nabla_{Y}(X)-[X,Y],
\end{equation*}
which is a skew-symmetric tensor with values in vector fields, i.e. $T\in\mathcal{D}(M)\otimes\Omega^{2}(M)$. In coordinates, it has the form
\begin{equation*}
T=\sum\limits_{i,j,k}(\Gamma_{ij}^{k}-\Gamma_{ji}^{k})\partial_{k}\otimes d_{i}\wedge d_{j}.
\end{equation*}
The connection is called \textit{torsion-free}, if $T=0$, i.e. $\Gamma_{ij}^{k}=\Gamma_{ji}^{k}$.

The \textit{curvature tensor} $C$ of a connection $\nabla$ is
\begin{equation*}
C\in\mathrm{End}(\mathcal{D}(M))\otimes\Omega^{2}(M),\quad C(X,Y)(Z)=[\nabla_{X},\nabla_{Y}]Z-\nabla_{[X,Y]}Z,
\end{equation*}
where $C(X,Y)\in\mathrm{End}(\mathcal{D}(M))$. In coordinates it has the form
\begin{equation*}
C=\sum\limits_{i,j,k,l}C^{i}_{jkl}\partial_{i}\otimes d_{j}\otimes d_{k}\wedge d_{l},
\end{equation*}
where coefficients $C^{k}_{lij}$ are related to Christoffel symbols by the following way:
\begin{equation*}
C_{lij}^{k}=\frac{\partial\Gamma_{lj}^{i}}{\partial x_{k}}-\frac{\partial\Gamma_{kj}^{i}}{\partial x_{l}}+\sum\limits_{m}(\Gamma_{lj}^{m}\Gamma_{km}^{i}-\Gamma_{kj}^{m}\Gamma_{lm}^{i}).
\end{equation*}
The torsion-free connection is said to be \textit{flat}, if $C=0$.

Let $(M,g)$ be a pseudo-Riemannian manifold with a pseudo-metric tensor $g$. Then, there exists a unique torsion-free connection, called \textit{Levi-Civita connection}, such that
\begin{equation*}
g(\nabla_{X}Y,Z)+g(Y,\nabla_{X}Z)=X(g(Y,Z)),\quad X,Y,Z\in\mathcal{D}(M).
\end{equation*}
This relation means that $\nabla_{X}(g)=0$ for all vector fields $X$. Christoffel symbols are related to metric $g$ as follows:
\begin{equation*}
\Gamma_{ij}^{k}=\frac{1}{2}\sum\limits_{l}g^{kl}\left(\frac{\partial g_{il}}{\partial x_{j}}+\frac{\partial g_{jl}}{\partial x_{i}}-\frac{\partial g_{ij}}{\partial x_{l}}\right),
\end{equation*}
where $g_{ij}=g(\partial_{i},\partial_{j})$ and $\|g^{ij}\|=\|g_{ij}\|^{-1}$.

Let $\mathcal{T}_{p}^{q}(M)=(\mathcal{D}(M))^{\otimes p}\otimes (\Omega^{1}(M))^{\otimes q}$ be the module of $p$-contravariant and $q$-covariant tensors on the manifold $M$ and let
\begin{equation*}
\mathcal{T}(M)=\oplus_{p,q}\mathcal{T}_{p}^{q}(M)
\end{equation*}
be the bigraded tensor algebra. Then, any affine connection $\nabla$ on the manifold $M$ defines a derivation $d_{\nabla}$ of degree $(1,1)$ in this algebra by the following way. On functions its action is $d_{\nabla}(f)=df$. Define this derivation on vector fields:
\begin{equation*}
d_{\nabla}\colon\mathcal{D}(M)\to\mathcal{D}(M)\otimes\Omega^{1}(M),\quad \langle d_{\nabla}(X),Y\rangle=\nabla_{Y}(X).
\end{equation*}
In coordinates we have
\begin{equation*}
d_{\nabla}(\partial_{i})=\sum\limits_{j,k}\Gamma^{k}_{ij}\partial_{k}\otimes d_{j}.
\end{equation*}
Then, we define this derivation on 1-forms:
\begin{equation*}
d_{\nabla}\colon\Omega^{1}(M)\to\Omega^{1}(M)\otimes\Omega^{1}(M),\quad d_{\nabla}(\omega)(Y,X)=X(\omega(Y))-\omega(\nabla_{X}(Y)).
\end{equation*}
In coordinates we have
\begin{equation*}
d_{\nabla}(d_{k})=-\sum\limits_{i,j}\Gamma^{k}_{ij}d_{j}\otimes d_{i}
\end{equation*}
The action of $d_{\nabla}$ on  higher order tensors is expanded by means of the Leibnitz rule:
\begin{equation*}
d_{\nabla}(\theta_{1}\otimes\theta_{2})=d_{\nabla}(\theta_{1})\otimes\theta_{2}+\theta_{1}\otimes d_{\nabla}(\theta_{2}).
\end{equation*}
We will use these constructions to get invariant symmetric tensors that will provide us with affine invariants on a plane.
\subsection{Symmetric Tensors}
Let $\Sigma^{k}(M)\subset(\Omega^{1}(M))^{\otimes k}$ be the module of symmetric tensors. Then,
\begin{equation*}
\Sigma^{\ast}(M)=\oplus_{k\ge0}\Sigma^{k}(M)
\end{equation*}
is a commutative algebra with the symmetric product. The derivation $d_{\nabla}$ defines a derivation of degree 1 in this algebra
\begin{equation*}
d_{\nabla}^{s}\colon\Sigma^{\ast}(M)\to\Sigma^{\ast+1}(M),
\end{equation*}
where
\begin{equation*}
d_{\nabla}^{s}\colon\Sigma^{k}(M)\stackrel{d_{\nabla}}{\longrightarrow}\Sigma^{k}(M)\otimes\Omega^{1}(M)\stackrel{\mathrm{Sym}}{\longrightarrow}\Sigma^{k+1}(M).
\end{equation*}
The derivation $\Sigma^{k}(M)$ allows to define higher order differentials $\theta_{k}(f)$ of functions $f\in C^{\infty}(M)$:
\begin{equation}
\label{sym-tens}
\Sigma^{k}(M)\ni\theta_{k}(f)=(d_{\nabla}^{s})^{k}(f)
\end{equation}
\begin{example}
Consider torsion-free connection $\nabla$. Then, we have
\begin{equation*}
\theta_{1}(f)=df=\sum\limits_{k}\partial_{k}(f)d_{k},
\end{equation*}
\begin{equation*}
\theta_{2}(f)=\sum\limits_{i,j}\left(\partial_{ij}(f)-\sum\limits_{k}\Gamma_{ij}^{k}\partial_{k}(f)\right)d_{i}\cdot d_{j}.
\end{equation*}
\end{example}
\subsection{Affine Invariants}
Let us consider affine invariants of the plane. The affine Lie algebra
\begin{equation*}
\mathfrak{aff}_{2}=\langle\partial_{x},\partial_{y},x\partial_{x},x\partial_{y},y\partial_{x},y\partial_{y}\rangle
\end{equation*}
acts transitively on $\mathbb{R}^{2}$, and therefore $\mathbf{J}^{k}/\mathfrak{aff}_{2}=\mathbf{J}^{k}_{0}/\mathfrak{gl}_{2}$, where
\begin{equation*}
\mathfrak{gl}_{2}=\langle x\partial_{x},x\partial_{y},y\partial_{x},y\partial_{y}\rangle.
\end{equation*}
The group of affine transformations preserves the trivial connection $\nabla^{s}$, therefore due to construction \eqref{sym-tens} symmetric tensors
\begin{equation*}
\Theta_{k}=\sum\limits_{i=0}^{k}u_{i,k-i}\frac{dx^{i}}{i!}\frac{dy^{k-i}}{(k-i)!}
\end{equation*}
are invariants of affine transformations.

Similar to Sect.\ref{sec:2}, we construct an invariant frame $\nabla_{1},\nabla_{2}$
\begin{equation*}
\nabla_{i}=A_{i}\frac{d}{dx}+B_{i}\frac{d}{dy},
\end{equation*}
such that
\begin{equation*}2\nabla_{1}\rfloor\Theta_{2}=\Theta_{1},\quad\Theta_{2}(\nabla_{1},\nabla_{2})=0,\quad\Theta_{2}(\nabla_{1},\nabla_{1})=\Theta_{2}(\nabla_{2},\nabla_{2}).
\end{equation*}
Then, we get
\begin{eqnarray*}
  \nabla_{1} &=&\frac{u_{02}u_{10}-u_{11}u_{01}}{u_{20}u_{02}-u_{11}^{2}}\frac{d}{dx}+\frac{u_{20}u_{01}-u_{11}u_{10}}{u_{20}u_{02}-u_{11}^{2}}\frac{d}{dy}, \\
  \nabla_{2} &=&\frac{1}{\sqrt{u_{20}u_{02}-u_{11}^{2}}}\left(-u_{01}\frac{d}{dx}+u_{10}\frac{d}{dy}\right),
\end{eqnarray*}
Note that the function $I_{0}=\Theta_{0}=u_{00}$ is an affine invariant of order zero, and therefore the function
\begin{equation*}
I_{2}=\nabla_{1}(I_{0})=\Theta_{1}(\nabla_{1})=2\Theta_{2}(\nabla_{1},\nabla_{1})=\|\nabla_{1}\|^{2}=\frac{u_{01}^{2}u_{20}-2u_{10}u_{01}u_{11}+u_{10}^{2}u_{02}}{u_{20}u_{02}-u_{11}^{2}}
\end{equation*}
is a second order differential affine invariant.

The dual coframe $\langle\omega_{1},\omega_{2}\rangle$ consists of horizontal 1-forms, such that $\omega_{i}(\nabla_{j})=\delta_{ij}$, and has the form
\begin{eqnarray*}
 \omega_{1} &=&\frac{1}{I_{2}}(u_{10}dx+u_{01}dy), \\
  \omega_{2} &=&\frac{1}{I_{2}\sqrt{u_{20}u_{02}-u_{11}^{2}}}\left((u_{11}u_{10}-u_{01}u_{20})dx+(u_{10}u_{02}-u_{11}u_{01})dy\right),
\end{eqnarray*}
and we also get an affine invariant volume form
\begin{equation*}
\omega_{1}\wedge\omega_{2}=\frac{\sqrt{u_{20}u_{02}-u_{11}^{2}}}{I_{2}}dx\wedge dy.
\end{equation*}

Summarizing above discussion, we observe that any regular function $f$ defines the following geometric structures associated with the affine geometry on $\mathbb{R}^{2}$
\begin{itemize}
\item pseudo-Riemannian structure $\Theta_{2}(f)$, that gives all Riemannian invariants \cite{LychYum1},
\item symplectic structure $(\omega_{1}\wedge\omega_{2})(f)$,
\item cubic form $\Theta_{3}(f)$ and Wagner connection \cite{LychYum2},
\end{itemize}
and others.

Writing down symmetric tensors $\Theta_{k}$ in terms of invariant coframe, we get
\begin{equation*}
\Theta_{k}=\sum\limits_{i=0}^{k}I_{i,k-i}\frac{\omega_{1}^{i}}{i!}\frac{\omega_{2}^{k-i}}{(k-i)!},
\end{equation*}
which gives us rational affine invariants (perhaps one should take squares to get rid of square roots) $I_{0}=u_{00}$,
\begin{equation}
I_{2}=\frac{u_{01}^{2}u_{20}-2u_{10}u_{01}u_{11}+u_{10}^{2}u_{02}}{u_{20}u_{02}-u_{11}^{2}},
\end{equation}
and $I_{i,k-i}$.

Since $\dim\mathbf{J}^{k}_{0}=\binom{k+2}{2}$ and $\dim(\mathfrak{gl}_{2})=4$ we observe that functions $I_{0},I_{2},I_{i,k-i}, 3\le i \le k$ generate the field of rational affine differential invariants of order $k$.
\subsection{Invariants of Algebraic Curves}
A plane algebraic curve is given by equation
\begin{equation*}
P_{k}(x,y)=0,
\end{equation*}
where $P_{k}(x,y)$ is an irreducible polynomial of degree $k$, which is defined up to a multiplier $P_{k}\mapsto\lambda P_{k}$, $\lambda\ne0$. This action is generated by an infinitely prolonged vector field $u_{00}\partial_{u_{00}}$:
\begin{equation*}
\gamma=\sum\limits_{ij}u_{ij}\frac{\partial}{\partial u_{ij}}.
\end{equation*}
An invariant $I$ is said to be of weight $w(I)$, if and only if
\begin{equation*}
\gamma(I)=w(I)I.
\end{equation*}
Affine invariants of zero weight are affine invariants of algebraic plane curves. Since $w(I_{0})=w(I_{2})=w(I_{i,j})=1$, one can choose
\begin{equation*}
\mathfrak{a}_{2}=\frac{I_{2}}{I_{0}},\quad\mathfrak{a}_{ij}=\frac{I_{ij}}{I_{0}}
\end{equation*}
as a generating set of rational affine invariants of algebraic plane curves.
\begin{remark}
An algebraic plane curve is defined by its $k$-th jet at the point $\mathbf{0}$, and therefore values
\begin{equation*}
\mathfrak{a}_{2}(P_{k})(0),\quad\mathfrak{a}_{ij}(P_{k})(0)
\end{equation*}
define the curve (completely over $\mathbb{C}$ and up to $\pm$ over $\mathbb{R}$).
\end{remark}
To find rational invariants (without square roots of the Hessian) we will use the coframe given by total differentials of invariants $I_{0}=u_{00}$  and $I_{2}=(u_{01}^{2}u_{20}-2u_{10}u_{01}u_{11}+u_{10}^{2}u_{02})(u_{20}u_{02}-u_{11}^{2})^{-1}$:
\begin{eqnarray*}
 \omega_{1} &=&\widehat{d}u_{00}=\Theta_{1}, \\
  \omega_{2} &=&\widehat{d}I_{2},
\end{eqnarray*}
and the Tresse frame as follows:
\begin{eqnarray*}
 \tau_{1} &=&A_{11}\frac{d}{dx}+A_{12}\frac{d}{dy}, \\
  \tau_{2} &=&A_{21}\frac{d}{dx}+A_{22}\frac{d}{dy},
\end{eqnarray*}
where
\begin{equation*}
\begin{pmatrix}
A_{11} & A_{12}\\
A_{21} & A_{22}
\end{pmatrix}=\begin{pmatrix}
u_{10} & \frac{dI_{2}}{dx}\\
u_{01} & \frac{dI_{2}}{dy}
\end{pmatrix}^{-1}.
\end{equation*}
Expressing the original coframe $\langle dx,dy\rangle$, we get
\begin{equation*}
\begin{pmatrix}
dx\\
dy
\end{pmatrix}=\begin{pmatrix}
u_{10} & u_{01}\\
\frac{dI_{2}}{dx} & \frac{dI_{2}}{dy}
\end{pmatrix}^{-1}\begin{pmatrix}
\omega_{1}\\
\omega_{2}
\end{pmatrix}.
\end{equation*}
Again, expression for symmetric tensors $\Theta_{k}$ in terms of the Tresse coframe
\begin{equation}
\label{inv-tens}
\Theta_{k}=\sum\limits_{i=0}^{k}I_{i,k-i}\frac{\omega_{1}^{i}}{i!}\frac{\omega_{2}^{k-i}}{(k-i)!},
\end{equation}
gives us affine invariants $I_{i,k-i}$ of the weight $(1-k)$, and we get
\begin{theorem}
Rational affine differential invariants are rational functions of invariants $I_{ij}$ given by \eqref{inv-tens}.
\end{theorem}
For algebraic curves, we have
\begin{theorem}
Rational affine differential invariants of algebraic curves are rational functions of invariants $I_{ij}I_{0}^{i+j-1}$.
\end{theorem}
\section{Invariants of Ternary Forms}
\label{sec:5}
In this section, we discuss the $\mathrm{SL}_{3}(\mathbb{C})$-classification problem for ternary forms of an arbitrary degree $n$, similar to the case of binary forms considered in Sect. \ref{sec:2}.

Ternary forms of degree $n$ are homogeneous polynomials on $\mathbb{C}^{3}$ of the form
\begin{equation}
\label{tern-form}
\mathcal{T}_{n}\ni\phi_{b}=\sum\limits_{i+j+k=n}b_{i,j,k}\frac{x^{i}}{i!}\frac{y^{j}}{j!}\frac{z^{k}}{k!}.
\end{equation}
The action of the Lie group
\begin{equation*}
\mathrm{SL}_{3}(\mathbb{C})=\left\{A\in\mathrm{Mat}_{3\times 3}(\mathbb{C})\mid \det(A)=1\right\}
\end{equation*}
on $\mathcal{T}_{n}$ is defined by the following way:
\begin{equation}
\label{SL3-act}
A\colon\mathcal{T}_{n}\ni\phi_{b}\mapsto A\phi_{b}=\phi_{b}\circ A^{-1}\in\mathcal{T}_{n}.
\end{equation}
The corresponding Lie algebra $\mathfrak{sl}_{3}$ consists of vector fields:
\begin{equation*}
X_{1}=x\partial_{x}-y\partial_{y},\quad X_{2}=x\partial_{x}-z\partial_{z},\quad X_{3}=y\partial_{x},\quad X_{4}=z\partial_{x},
\end{equation*}
\begin{equation*}
X_{5}=x\partial_{y},\quad X_{6}=z\partial_{y},\quad X_{7}=x\partial_{z},\quad X_{8}=y\partial_{z}.
\end{equation*}
Similar to the case of binary forms, we consider \eqref{tern-form} as smooth solutions to the Euler equation:
\begin{equation}
\label{Euler-ternary}
xf_{x}+yf_{y}+zf_{z}=nf.
\end{equation}
Equation \eqref{Euler-ternary} defines a smooth manifold in the space of $1$-jets of functions on $\mathbb{C}^{3}$:
\begin{equation*}
\mathcal{E}_{1}=\left\{xu_{100}+yu_{010}+zu_{001}=nu_{000}\right\}\subset\mathbf{J}^{1}.
\end{equation*}
As in the previous sections, we will use the notation $\mathcal{E}_{k}$ for the collection of all prolongations of \eqref{Euler-ternary} to the space $\mathbf{J}^{k}$ up to order $k$.

The action $A\colon\mathbb{C}^{3}\to\mathbb{C}^{3}$ of the group $\mathrm{SL}_{3}$ can be prolonged to $\mathbf{J}^{k}$ by the natural way
\begin{equation*}
A^{(k)}\colon\mathbf{J}^{k}\to\mathbf{J}^{k},\quad A^{(k)}\left([f]_{p}^{k}\right)=[Af]_{Ap}^{k}.
\end{equation*}

A rational function $I\in C^{\infty}(\mathcal{E}_{k})$ is said to be a \textit{differential $\mathrm{SL}_{3}$-invariant of order $k$}, if $I\circ A^{(k)}=I$, for all $A\in\mathrm{SL}_{3}(\mathbb{C})$.

Using the results of Sect. \ref{sec:4} we define $\mathrm{SL}_{3}(\mathbb{C})$-invariant symmetric tensors:
\begin{equation}
\label{inv-tens-tern}
\Theta_{m}=\sum\limits_{i+j+k=m}u_{ijk}\frac{dx^{i}}{i!}\frac{dy^{j}}{j!}\frac{dz^{k}}{k!}.
\end{equation}
To construct an invariant coframe we will need an inverse of $\Theta_{2}$:
\begin{equation*}
\begin{split}
\Theta_{2}^{-1}&=\frac{2}{A}((u_{002}u_{020}-u_{011}^{2})\partial_{x}\partial_{x}-2(u_{002}u_{110}-u_{011}u_{101})\partial_{x}\partial_{y}+{}\\&+2(u_{011}u_{110}-u_{020}u_{101})\partial_{x}\partial_{z}-2(u_{011}u_{200}-u_{101}u_{110})\partial_{y}\partial_{z}+{}\\&+(u_{002}u_{200}-u_{101}^{2})\partial_{y}\partial_{y}+(u_{020}u_{200}-u_{110}^{2})\partial_{z}\partial_{z}),
\end{split}
\end{equation*}
where
\begin{equation*}
A=u_{002}u_{020}u_{200}-u_{002}u_{110}^{2}-u_{011}^{2}u_{200}+2u_{011}u_{101}u_{110}-u_{020}u_{101}^{2}
\end{equation*}
is a differential $\mathrm{SL}_{3}(\mathbb{C})$-invariant of order 2.

As the first invariant form $\omega_{1}$, we take
\begin{equation*}
\omega_{1}=\Theta_{1}=u_{100}dx+u_{010}dy+u_{001}dz.
\end{equation*}
The second invariant form will be the total differential of the invariant $A$
\begin{equation*}
\omega_{2}=\frac{dA}{dx}dx+\frac{dA}{dy}dy+\frac{dA}{dz}dz=A_{1}dx+A_{2}dy+A_{3}dz,
\end{equation*}
where
\begin{equation*}
\begin{split}
A_{1}&=u_{002}u_{020}u_{300}-2u_{002}u_{110}u_{210}+u_{002}u_{120}u_{200}-u_{011}^{2}u_{300}+{}\\&+2u_{011}u_{101}u_{210}+2u_{011}u_{110}u_{201}-2u_{011}u_{111}u_{200}-2u_{020}u_{101}u_{201}+{}\\&+u_{020}u_{102}u_{200}-u_{101}^{2}u_{120}+2u_{101}u_{110}u_{111}-u_{102}u_{110}^{2}
\end{split}
\end{equation*}
\begin{equation*}
\begin{split}
A_{2}&=u_{002}u_{020}u_{210}+u_{002}u_{030}u_{200}-2u_{002}u_{1 1 0}u_{1 2 0}-u_{0 1 1} ^{2}u_{2 1 0}-2u_{0 1 1}u_{0 2 1}u_{2 0 0}+{}\\&+2u_{0 1 1}u_{1 0 1}u_{1 2 0}+2u_{0 1 1}u_{1 1 0}u_{1 1 1}+u_{0 1 2}u_{0 2 0}u_{2 0 0}-u_{0 1 2}u_{1 1 0} ^{2}-2u_{0 2 0}u_{1 0 1}u_{1 1 1}+{}\\&+2u_{0 2 1}u_{1 0 1}u_{1 1 0}-u_{0 3 0}u_{1 0 1} ^{2}
\end{split}
\end{equation*}
\begin{equation*}
\begin{split}
A_{3}&=u_{0 0 2}u_{0 2 0}u_{2 0 1}+u_{0 0 2}u_{0 2 1}u_{2 0 0}-2u_{0 0 2}u_{1 1 0}u_{1 1 1}+u_{0 0 3}u_{0 2 0}u_{2 0 0}-{}\\&-u_{0 0 3}u_{1 1 0}^{2}-u_{0 1 1}^{2}u_{2 0 1}-2u_{0 1 1}u_{0 1 2}u_{2 0 0}+2u_{0 1 1}u_{1 0 1}u_{1 1 1}+2u_{0 1 1}u_{1 0 2}u_{1 1 0}+{}\\&+2u_{0 1 2}u_{1 0 1}u_{1 1 0}-2u_{0 2 0}u_{1 0 1}u_{1 0 2}-u_{0 2 1}u_{1 0 1}^{2}.
\end{split}
\end{equation*}
The third invariant form $\omega_{3}=F_{1}dx+F_{2}dy+F_{3}dz$ is found from the conditions of orthogonality to $\omega_{2}$ and $\Theta_{1}$ in the sense of $\Theta_{2}$:
\begin{equation*}
\Theta_{2}^{-1}(\omega_{2},\omega_{3})=0,\quad\Theta_{2}^{-1}(\Theta_{1},\omega_{3})=0,
\end{equation*}
which define the form $\omega_{3}$ up to a multiplier:
\begin{eqnarray*}
F_{1}&=&F_{3}\frac{(u _{0 0 1}u _{1 1 0}-u _{0 1 0}u _{1 0 1})A_{1}+(-u _{0 0 1}u _{2 0 0}+u _{1 0 0}u _{1 0 1})A_{2}+(u _{0 1 0}u _{2 0 0}-u _{1 0 0}u _{1 1 0})A_{3}}{(u _{0 0 1}u _{0 1 1}-u _{0 0 2}u _{0 1 0})A_{1}+(-u _{0 0 1}u _{1 0 1}+u _{0 0 2}u _{1 0 0})A_{2}+(u _{0 1 0}u _{1 0 1}-u _{0 1 1}u _{1 0 0})A_{3}},\\
F_{2}&=&F_{3}\frac{(u_{0 0 1} u_{0 2 0}-u_{0 1 0} u_{0 1 1})A_{1}+(-u_{0 0 1} u_{1 1 0}+u_{0 1 1} u_{1 0 0})A_{2}+(u_{0 1 0} u_{1 1 0}-u_{0 2 0} u_{1 0 0})A_{3}}{(u_{0 0 1} u_{0 1 1}-u_{0 0 2} u_{0 1 0})A_{1}+(-u_{0 0 1} u_{1 0 1}+u_{0 0 2} u_{1 0 0})A_{2}+(u_{0 1 0} u_{1 0 1}-u_{0 1 1} u_{1 0 0})A_{3}}.
\end{eqnarray*}
We put $F_{3}$ equal to the denominator in the above expressions:
\begin{equation*}
F_{3}=(u_{0 0 1} u_{0 1 1}-u_{0 0 2} u_{0 1 0})A_{1}+(-u_{0 0 1} u_{1 0 1}+u_{0 0 2} u_{1 0 0})A_{2}+(u_{0 1 0} u_{1 0 1}-u_{0 1 1} u_{1 0 0})A_{3}.
\end{equation*}
One can check that in this case the form $\omega_{3}$ will be invariant.

Now that we have constructed an invariant coframe $\langle\omega_{1},\omega_{2},\omega_{3}\rangle$, we are able to construct an invariant frame $\langle\nabla_{1},\nabla_{2},\nabla_{3}\rangle$ dual to $\langle\omega_{1},\omega_{2},\omega_{3}\rangle$:
\begin{equation*}
\omega_{i}(\nabla_{j})=\delta_{ij}.
\end{equation*}
And finally we are able to express the original coframe $\langle dx,dy,dz\rangle$ in terms of an invariant one:
\begin{equation*}
\begin{pmatrix}
dx \\ dy\\ dz
\end{pmatrix}=
\begin{pmatrix}
u_{100} & u_{010} & u_{001}\\
\frac{dA}{dx} & \frac{dA}{dy} & \frac{dA}{dz}\\
F_{1} & F_{2} & F_{3}
\end{pmatrix}^{-1}\begin{pmatrix}
\omega_{1}\\ \omega_{2} \\ \omega_{3}
\end{pmatrix}.
\end{equation*}

Therefore tensors \eqref{inv-tens-tern} are written by the following way:
\begin{equation*}
\Theta_{m}=\sum\limits_{i+j+k=m}I_{ijk}\frac{\omega_{1}^{i}}{i!}\frac{\omega_{2}^{j}}{j!}\frac{\omega_{3}^{k}}{k!}.
\end{equation*}
\begin{theorem}
Functions $I_{ijk}$ are $\mathrm{SL}_{3}$-invariants of order $(i+j+k)$, and any rational differential invariant is a rational function of them.
\end{theorem}

However, explicit expressions for invariants $I_{i,j,k}$ look bulky and straightforward computations work slowly in the case of ternary forms. To this reason, to find a generating set of invariants, we will use the Lie-Tresse theorem. Namely, we take five third-order independent invariants
\begin{equation}
\label{invs-tern}
J_{1}=u_{00},\quad J_{2}=A,\quad J_{3}=\nabla_{1}(J_{2}),\quad J_{4}=\nabla_{2}(J_{2}),\quad J_{5}=\nabla_{3}(J_{2}).
\end{equation}
Since $\dim\mathcal{E}_{3}=13$, $\dim\mathfrak{sl}_{3}=8$, then we need five differential invariants to separate regular orbits. According to the global Lie-Tresse theorem, all other rational differential invariants can be found from \eqref{invs-tern} by applying invariant derivations $\nabla_{i}$.
\begin{theorem}
The field of rational $\mathfrak{sl}_{3}$-invariants is generated by \eqref{invs-tern} and invariant derivations $\nabla_{i}$. They separate regular orbits.
\end{theorem}
If we restrict \eqref{invs-tern} to the ternary form of degree $n$, we will get five functions on a three-dimensional space, therefore, there are 2 relations between them:
\begin{equation}
\label{quo-tern}
F_{1}(J_{1}^{\phi},J_{2}^{\phi},J_{3}^{\phi},J_{4}^{\phi},J_{5}^{\phi})=0,\quad F_{2}(J_{1}^{\phi},J_{2}^{\phi},J_{3}^{\phi},J_{4}^{\phi},J_{5}^{\phi})=0.
\end{equation}
To write out syzygies \eqref{quo-tern} explicitly, one can use the similar Maple code as we used in Sect. \ref{sec:2} for cubics.
\begin{theorem}
Let $\phi$ be a regular ternary form of degree $n$. Then, $\mathrm{SL}_{3}(\mathbb{C})$-orbit of $\phi$ consists of solutions to a quotient PDE
\begin{equation*}
F_{1}(J_{1},J_{2},J_{3},J_{4},J_{5})=0,\quad F_{2}(J_{1},J_{2},J_{3},J_{4},J_{5})=0.
\end{equation*}
together with $\mathcal{E}_{n}$.
\end{theorem}
\section*{Acknowledgements}
This work was partially supported by the Foundation for the Advancement of Theoretical Physics and Mathematics ``BASIS'' (project 19-7-1-13-3).

\end{document}